\documentclass[11pt]{article}
\usepackage[margin=1in]{geometry}
\usepackage{amsmath,amssymb,booktabs,graphicx,hyperref,xcolor}
\usepackage[round,authoryear]{natbib}
\usepackage{subcaption}
\usepackage{needspace}

\title{When Language Models Misinterpret Uncertainty in Earth System Mechanism Attribution}
\author{Pouria Behnoudfar,\qquad  Nan Chen\\
 Department of Mathematics, University of Wisconsin-Madison}

\date{}

\begin{document}
	\maketitle
	

\section*{Key Points}

\begin{itemize}
    \item Providing uncertainty statistics does not ensure that language models correctly judge competing physical mechanisms.

    \item Numerical separation improved one language model in idealized tests but not clearly in coupled climate tests.

    \item Combining complete uncertainty information with an explicit interpretation reduced false ambiguity reports in coupled climate tests.
\end{itemize}

\section*{Plain Language Summary}

An unusual ocean or climate event may have several physical explanations. Scientists use observations and mathematical models to estimate the role of each process, but those estimates can be uncertain. When language models help describe scientific analyses, do they convey that uncertainty? We investigated this using two language models, a simple mathematical system, and a model of El Ni\~{n}o in the Pacific. The models received different kinds of evidence, from a single proposed mechanism to numerical uncertainty measures and a full assessment accompanied by its interpretation. More numbers alone did not consistently improve their answers. With the full assessment and interpretation, both models more reliably reported when the evidence could not clearly distinguish competing processes. This also reduced cases incorrectly described as ambiguous. Applied to historical Pacific events, the procedure often reproduced model-derived diagnoses, without establishing that these were the actual physical causes. The findings highlight why uncertainty assessment and its explanation should be evaluated separately in AI-assisted Earth system research.

\begin{abstract}

Mechanistic attribution in Earth system science requires quantifying uncertainty and communicating its implications. Here, we develop an uncertainty-aware mechanistic attribution pipeline that constructs probabilistic fingerprints from ensemble estimates, including covariance-aware separation of competing mechanisms, and an evidence ladder that systematically varies the information provided to language models (LLMs). We evaluate this framework in an idealized stochastic system and a coupled El Niño model using controlled interventions and ensemble data assimilation. In the idealized system, numerical separation increases GPT-4o's intervention-label agreement from 0.80 to 0.98, but not GPT-4o-mini's. In the coupled system, both LLMs identify seven of eight fingerprint-ambiguous events when given full fingerprints with verbal interpretations; GPT-4o's false ambiguity reports decrease from 10/82 to 1/82, without consistent gains in intervention-label agreement. For eight historical ENSO events, the uncertainty-aware condition agrees with the model-inferred mechanism labels in seven cases. These findings distinguish uncertainty assessment from its communication.

\end{abstract}

	\section{Introduction}
	\label{sec:intro}

Earth system anomalies may arise through different physical pathways, making it difficult to identify their causes from observed changes alone. Marine heatwaves, for example, can involve air--sea exchange and ocean transport, while droughts can reflect atmospheric circulation and land--atmosphere feedbacks~\citep{holbrook2019global,stott2016,shepherd2018,hobday2016}. In the tropical Pacific, El Ni\~{n}o events can be influenced by thermocline feedback, zonal advection, and westerly wind bursts~\citep{vialard2025nino,capotondi2015,timmermann2018,chen2015wwb}. Identifying the leading mechanism is not the same as establishing that it can be distinguished from alternatives.

Making such attributions is challenging even when the observed anomaly is well characterized. The variables used to describe an event do not necessarily identify the physical processes responsible for it. Mechanism-related states, fluxes, forcings, and parameters often need to be estimated from observations and models. These estimates carry uncertainty arising from incomplete measurements, internal variability, uncertain parameters, and model assumptions. Different plausible reconstructions may therefore be consistent with the same observations while assigning different relative importance to the candidate mechanisms \citep{tarantola2005,ruiz2013,poulet2024slip}. Moreover, several mechanisms may genuinely contribute to an event, independently of the uncertainty in their estimated contributions~\citep{runge2019}. A single best estimate may identify the mechanism with the largest inferred contribution, but it does not establish how robust that ranking is or whether other mechanisms remain comparably plausible~\citep{hawkins2009}.

Uncertainty quantification (UQ) provides a way to assess the strength of such attributions. Depending on the problem, uncertainty in mechanism-related quantities can be characterized through ensemble prediction \citep{palmer2000,gneiting2005}, Bayesian inference \citep{tarantola2005}, data assimilation \citep{evensen1994,kalnay2003,vanleeuwen2009,law2015,behnoudfar2026rl,behnoudfar2026bridging}, or ensembles with uncertain parameters and forcings \citep{giorgini2025predicting, santos2026nonequilibrium,behnoudfar2026omega,behnoudfar2026split}. Rather than evaluating a mechanism diagnostic at a single estimate, one can examine its distribution across plausible realizations. The resulting statistics show not only which mechanism has the highest average score, but also how clearly that score is separated from those of the alternatives, given their uncertainty. Reporting only the leading mechanism can therefore convey more certainty than the evidence supports.

Communicating these distinctions presents a further challenge. Large language models (LLMs) are increasingly used to explain scientific findings about Earth system phenomena \citep{vaghefi2023,bi2024}. Their ability to produce coherent explanations \citep{bulian2023assessing,lu2025can}, however, does not guarantee that those explanations accurately reflect the strength of the underlying evidence \citep{ji2023,kadavath2022,xiong2024,vanderbles2019}. If an LLM receives only a point estimate of the mechanism contributions, it lacks the information needed to assess how sensitive that attribution is to uncertainty. Providing quantitative uncertainty information addresses this limitation, but does not necessarily ensure that the model interprets it correctly \citep{perez2025uncertainty}. For example, the spread of individual mechanism scores does not directly determine whether two mechanisms can be distinguished; their separation also depends on the score difference and the joint variability of the scores. Without considering this relationship, an LLM may identify the highest-ranking mechanism without adequately assessing whether the available evidence supports a definitive attribution.

For AI-assisted Earth system analysis, the risk is not only an incorrect mechanism label: a fluent explanation may describe a poorly resolved diagnosis as definitive, or incorrectly report ambiguity \citep{camps2026bridging}. We therefore study how numerical uncertainty should be presented so that a language model conveys what the mechanism diagnostics actually distinguish, rather than simply producing a plausible explanation~\citep{turpin2023,jacovi2020}.

Here, we investigate this question through an uncertainty-aware mechanistic attribution pipeline and an evidence ladder that systematically varies the information provided to an LLM. The pipeline summarizes uncertain mechanism diagnostics in probabilistic fingerprints, while the evidence ladder examines how different quantitative and verbal representations affect the resulting attribution. We study the framework in two complementary settings: an idealized stochastic system and an intermediate coupled model of El Ni\~{n}o variability. In controlled twin experiments, selected physical processes are amplified to define intervention-based reference labels, while comparisons with unperturbed simulations illustrate their effects on the observed variables. In the coupled-model experiments, we use ensemble data assimilation of SST observations to characterize uncertainty in mechanism-related physical fields. We then examine two language models under eight evidence conditions, evaluating both their agreement with the intervention-based labels and their ability to recognize when the probabilistic fingerprints do not clearly distinguish competing mechanisms. Finally, we apply the framework to historical El Ni\~{n}o and La Ni\~{n}a events using GODAS ocean reanalysis data. Without prescribed intervention labels for these events, we assess how faithfully the language models reproduce the attributions supported by the inferred fingerprints, rather than treating those attributions as independently verified physical explanations.


\section{Framework and Experimental Design}
\label{sec:framework}

\begin{figure}[htbp]
		\centering
		\includegraphics[width=\linewidth]{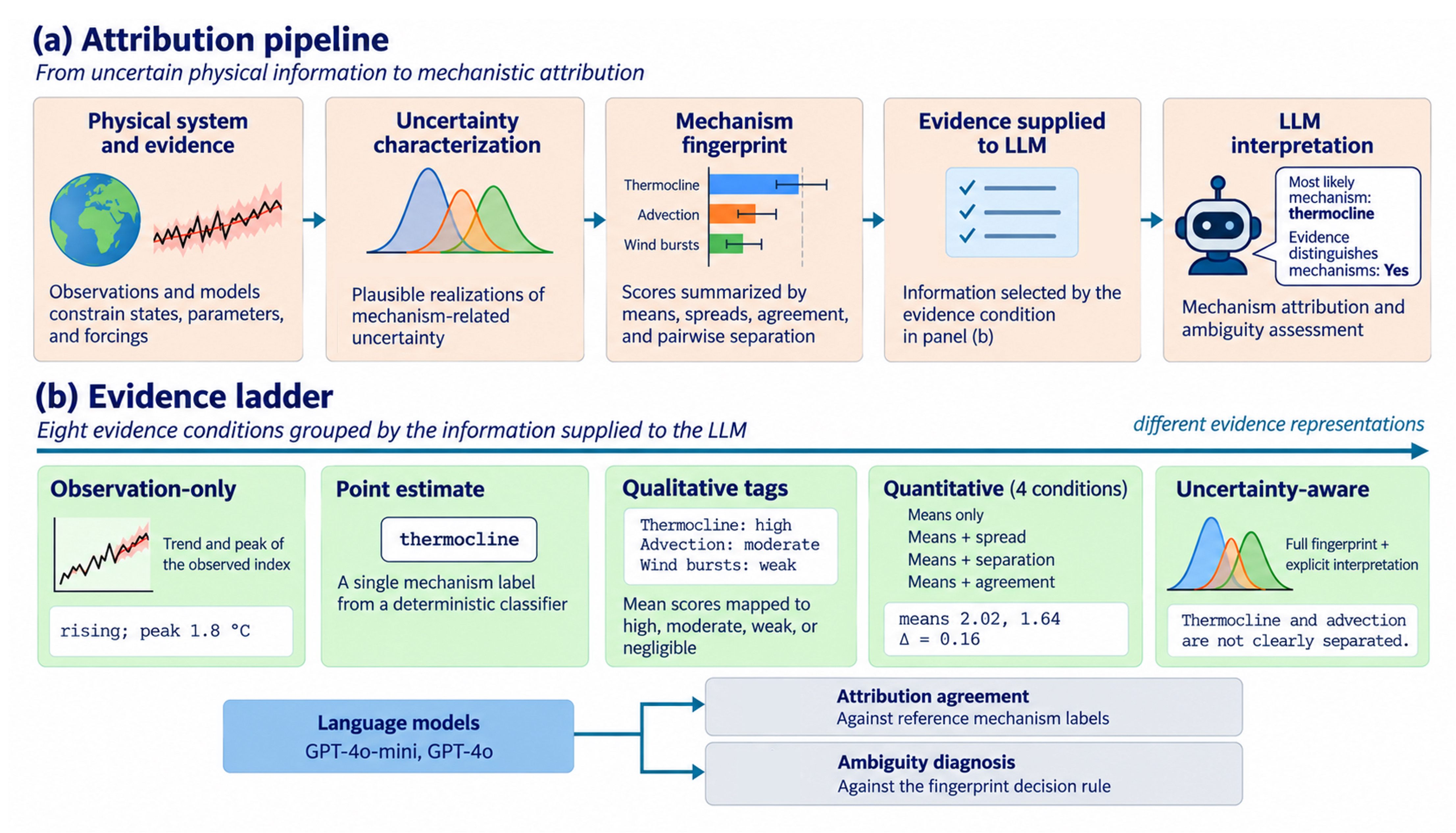}
		
\caption{
Overview of the framework and evaluation design. (a) Uncertainty in mechanism-related quantities is characterized and summarized in probabilistic fingerprints. An LLM interprets the evidence supplied from these fingerprints and produces an attribution. (b) Eight evidence conditions, grouped into five levels, examine how different representations of mechanism information affect LLM responses. The quantitative level includes a means-only condition and three variants that add one statistic at a time. Model settings and output format are held fixed across conditions.
}
		\label{fig:overview}
	\end{figure}

The framework couples an uncertainty-aware attribution pipeline to an evidence ladder for LLM interpretation (Figure~\ref{fig:overview}). The pipeline uses uncertain mechanism-related quantities to construct a fingerprint of relative support and distinguishability; the ladder varies the information presented to an LLM. Fingerprint labels and, in the controlled experiments, intervention-based reference labels are established independently of the LLM and need not coincide.

\subsection{Probabilistic mechanism fingerprints}

A mechanism fingerprint summarizes the evidence for different physical processes while retaining the uncertainty in their estimated activity. The starting point is an ensemble of plausible realizations of the mechanism-related quantities. These realizations may come from ensemble forecast, data assimilation, Bayesian inference, or other UQ methods, provided that their joint variability is represented.

For each candidate mechanism $k$, a diagnostic indicator $S_k$ is evaluated for every ensemble member $m=1,\ldots,M$. The indicator characterizes the activity or inferred influence of the corresponding physical process. Its value is standardized relative to a reference distribution obtained under background conditions:
\[
Z_k^{(m)}
=
\frac{S_k^{(m)}-\mu_k^{\rm clim}}
{\sigma_k^{\rm clim}},
\]
where $\mu_k^{\rm clim}$ and $\sigma_k^{\rm clim}$ are the reference mean and standard deviation (usually at the climatology). A positive score therefore indicates activity above the reference level, while a negative score indicates activity below it. In the experiments below, the indicators summarize mechanism-related activity over an event window rather than at a single time.

The fingerprint contains four summaries. (i) The mean score $\bar Z_k$ measures the average strength of each indicator relative to the reference conditions. (ii) Its ensemble spread $\sigma_{Z_k}$ characterizes uncertainty in that score. (iii) Ensemble agreement is the fraction of realizations in which mechanism $k$ receives the highest score. (iv) The pairwise separation compares the mean scores of two mechanisms relative to the uncertainty in their difference:
\[
\Delta_{ij}
=
\frac{\bar Z_i-\bar Z_j}
{\sqrt{
\sigma_{Z_i}^2+\sigma_{Z_j}^2
-2\,\mathrm{Cov}(Z_i,Z_j)
}}.
\]
The covariance is important because two mechanisms may respond together to the same uncertain physical state. A large $|\Delta_{ij}|$ indicates that their mean scores are well separated relative to the variability of their difference; a small value indicates that their relative ranking is less clearly resolved.

To obtain a fingerprint-based classification, we apply fixed decision rules to these scores. A mechanism is considered supported when $\bar Z_k>z_{\rm thr}$. When multiple mechanisms are supported, their separation determines whether the fingerprint favors a single leading mechanism or assigns a mixed label. If none is supported, the label is \emph{none}. Here, \emph{mixed} is an operational fingerprint label indicating that multiple mean-supported mechanisms cannot be clearly ranked under the adopted rule. It does not establish comparable causal contributions, and a mean-supported indicator may itself have substantial ensemble spread.

The support threshold is set using background simulations, while the separation threshold is chosen separately for the two numerical test systems. We use $z_{\rm thr}=1.0$ in both systems, with $\delta_{\rm thr}=2.5$ for the idealized model and $1.0$ for the ENSO model. The decision rules and threshold choices are described in Appendix~\ref{app:thresholds}.

\subsection{Test systems and event construction}

We examine the framework in two settings. The first is an idealized stochastic model used to establish a proof of concept. It combines three independent Ornstein-Uhlenbeck processes with different characteristic timescales that drive a common observed variable \citep{uhlenbeck1930}. By selectively amplifying the forcing of individual processes, we construct events with prescribed driver labels and examine how uncertainty affects their attribution. Appendix~\ref{app:concept} details how the mechanism indicators and the separation statistic label events in the conceptual model. 

The second setting uses a stochastic intermediate coupled El Ni\~no-Southern Oscillation (ENSO) model \citep{chen2023}. The model describes coupled ocean-atmosphere variability in the equatorial Pacific, with sea surface temperature (SST) represented on 56 longitudinal grid points and ocean dynamics represented through leading equatorial Kelvin and Rossby wave amplitudes. It reproduces aspects of observed ENSO variability, including the distributions and temporal autocorrelations of SST anomalies in different Ni\~{n}o regions (Figure~\ref{fig:stat}). We use this model to examine thermocline feedback, zonal advection, and wind-burst forcing, whose effects can overlap through the coupled dynamics.

In the ENSO experiments, we deliberately assimilate SST alone to study attribution under partial observability. Thermocline depth, zonal currents, and wind forcing are inferred through the stochastic coupled dynamics rather than directly observed in the assimilation. Although reanalysis products provide other variables, restricting the input to SST allows us to examine uncertainty in these model-inferred mechanism diagnostics.

\begin{figure}[htbp]
		\centering
		\includegraphics[width=1.0\textwidth]{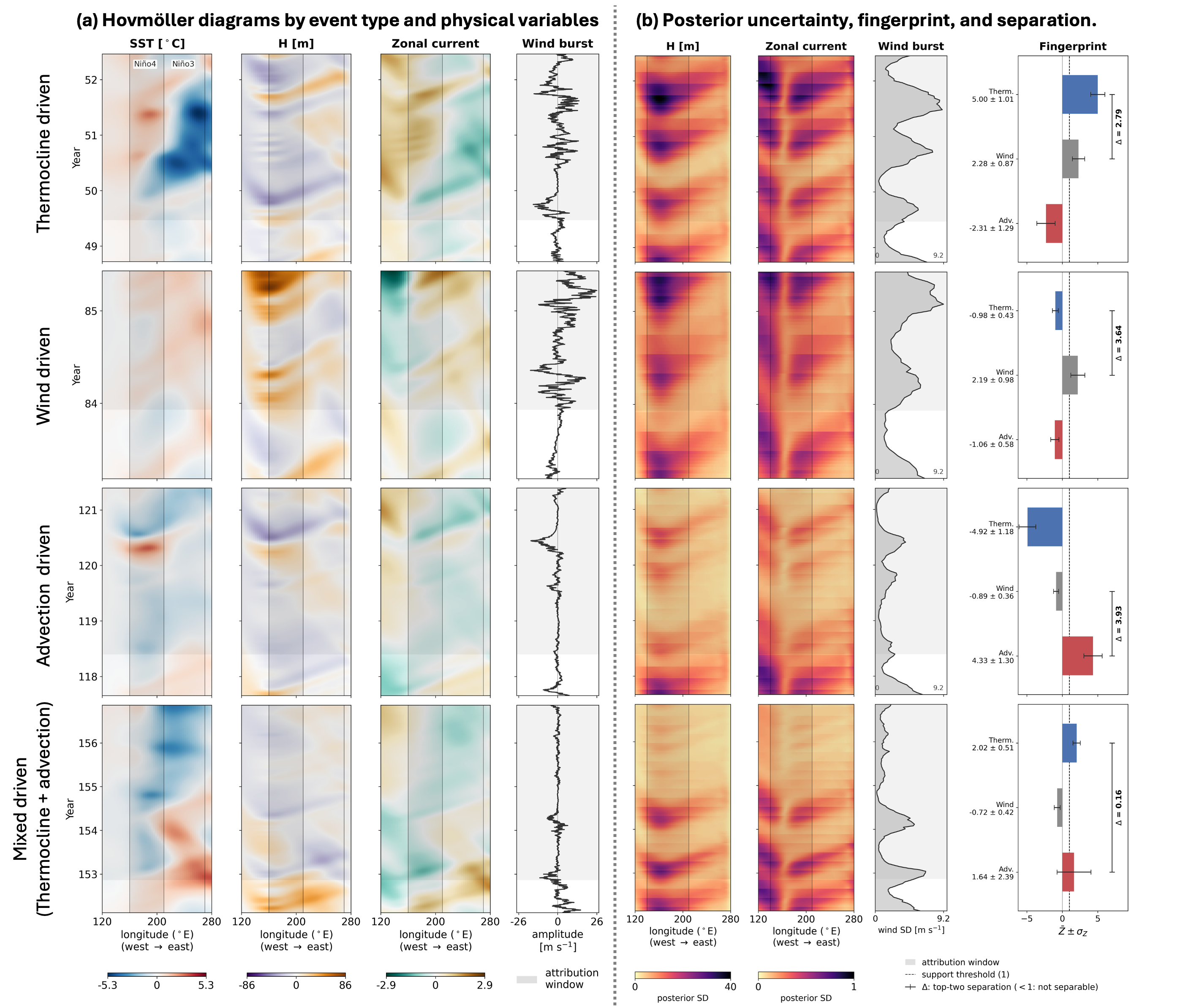}
		\caption{Controlled ENSO event library. (a) Physical fields from selected events in each event class whose fingerprints agree with the reference labels. (b) Posterior uncertainty of the four events  (rows: thermocline-driven, wind-burst-driven, advective-driven, mixed). The uncertainty is in physical units and on separate color scales. Columns: posterior SD of the thermocline depth $H$ and of the zonal current $U$ (Hovm\"oller diagrams), posterior SD of the wind-burst amplitude against time, and the fingerprint $\bar Z_k\pm\sigma_{Z_k}$ with the separation $\Delta$ of the two highest-ranked mechanisms. Color scales and axes are shared across columns and rows, respectively. The shaded band is the attribution window; the shaded longitude band is the Ni\~no-4 region, and the vertical lines bound Ni\~no-3. Dashed line: support threshold $z_{\rm thr}=1$.
  }
		\label{fig:saicm-diagnostics-app}
	\end{figure}

For the controlled experiments, each event is generated from an independent simulation. Following an initial background period, one or two selected physical processes are amplified over a prescribed event window. We compare perturbed and unperturbed simulations initialized from the same state and driven by the same realization of stochastic forcing. The reference label is determined by which processes are deliberately amplified, not by a threshold on the observed response. Paired simulations illustrate how these interventions affect the observed quantity (SST in the ENSO model), without establishing a unique causal decomposition of each event. The perturbations are intended to preserve the characteristic dynamics of the underlying system. Figure~\ref{fig:saicm-diagnostics-app}a illustrates selected events and their associated physical fields. Figure~\ref{fig:saicm-diagnostics-app}b shows the standard deviation of each variable in physical units together with the indicator scores and separation. Additional examples used to examine individual fingerprints are shown in Figure~\ref{fig:uq_blocka} with their corresponding uncertainty in Figure~\ref{fig:uq_blockc}.

We further apply the framework to historical El Ni\~{n}o and La Ni\~{n}a events using monthly GODAS SST for 1981--2020 \citep{behringer2004}. These SST fields are assimilated into the intermediate coupled model using the SST-only setup introduced above, yielding model-inferred mechanism diagnostics and their uncertainty. We analyze eight events selected to represent eastern- and central-Pacific El Ni\~{n}o and La Ni\~{n}a conditions. Mechanism fingerprints are computed over fixed 12-month windows, allowing the diagnostics to account for the evolution of the relevant processes around each event. The selection criteria, event windows, and score normalization are detailed in Appendix~\ref{app:godas_events}.
Unlike the controlled experiments, the historical events have no prescribed intervention labels. We therefore use the fingerprint-based attribution as the reference for evaluating LLM responses. This application assesses whether the language models faithfully communicate the conclusions supported by the model-inferred fingerprints, without treating those conclusions as independently verified physical attributions.

\subsection{Uncertainty quantification in the experiments}

For the ENSO experiments, we use a stochastic ensemble Kalman filter \citep{evensen1994,burgers1998} to estimate the model state and its uncertainty from SST input. The ensemble contains $M=150$ members. In the controlled experiments, SST is observed at the model's 56 longitudinal grid points, while the ocean current, thermocline depth, and wind forcing are inferred through the coupled dynamics. The reconstructed ensemble is then used to evaluate the mechanism indicators and construct their probabilistic fingerprints.

In the controlled experiments, the assimilating model uses smaller perturbation amplitudes than the reference simulations. This introduces a systematic difference between the reference and assimilating dynamics. Observations are assimilated monthly to obtain the ensemble used for attribution. The details of the filter, perturbation amplitudes, observation errors, and sensitivity to observation frequency are provided in Appendix~\ref{app:density}; the GODAS observation coverage and normalization are described in Appendix~\ref{app:godas_events}.

\subsection{Evidence ladder}

We use GPT-4o-mini, a smaller and more cost-efficient model, and GPT-4o, its more capable counterpart. Their comparison probes the dependence of attribution statements on the language model used. Each event is evaluated under eight evidence conditions, organized into five levels in Figure~\ref{fig:overview}b. The system instructions, generation settings, and structured output format are held fixed across conditions. Each event is submitted once under each condition through an independent request, allowing paired comparisons of the responses.

The first condition, \emph{observation-only}, provides the trend and peak of the observed scalar index (an SST index in the ENSO experiments), without mechanism diagnostics. The second, \emph{point estimate}, supplies only the mechanism selected by the deterministic classifier. The \emph{qualitative tags} condition describes the mean standardized scores as high, moderate, weak, or negligible. Four \emph{quantitative} conditions supply the numerical mean scores, either alone or with one additional statistic: ensemble spread, pairwise separation, or ensemble agreement. Finally, the \emph{uncertainty-aware} condition supplies the full fingerprint together with a short verbal statement of the attribution conclusion derived from it. This condition tests how faithfully the LLM conveys a supplied diagnosis; it does not isolate the effect of the verbal statement from that of the additional numerical fields.

The quantitative conditions retain some common ensemble-derived information, including relative support weights and, for the ENSO experiments, SST diversity diagnostics. Consequently, the comparisons that add one statistic at a time assess its contribution relative to this shared baseline. The qualitative and quantitative conditions also differ in some contextual information, so their comparison cannot be interpreted solely as an effect of replacing qualitative values with numerical ones. Full prompt specifications are given in Appendix~\ref{app:ladder-full}.

For each condition, the LLM is asked to identify the mechanism supported by the evidence and indicate whether the attribution is unresolved. We evaluate its agreement with the prescribed mechanism labels in the controlled experiments and separately assess whether it recognizes ambiguity according to the fingerprint rule. For the historical ENSO cases, where no intervention-based reference is available, we evaluate agreement with the model-inferred fingerprint labels.
	

\section{Results}
\label{sec:results}

Table~\ref{tab:ladder-main} summarizes two aspects of the LLM responses in the controlled experiments. Panel (a) reports agreement between the LLM attribution and the intervention-based reference label, defined by the processes deliberately perturbed in each experiment. Panel (b) examines whether the LLM recognizes events for which the probabilistic fingerprint does not clearly separate competing mechanisms. The measures need not agree: matching an intervention label is different from identifying ambiguity under the fingerprint rule. Because a fingerprint can itself differ from the intervention label, the agreement in panel (a) does not isolate LLM interpretation errors from errors in the preceding diagnostic step. The idealized events are the 60 events of the conceptual-model library with no gate, $Y_1$, $Y_2$, or both gated.

\begin{table}[htbp]
\centering
\small

\textit{(a) Agreement with mechanism reference labels}

\vspace{.2cm}
\begin{tabular}{@{}lcccc@{}}
\toprule
& \multicolumn{2}{c}{Idealized ($n=60$)}
& \multicolumn{2}{c}{Coupled ENSO ($n=90$)} \\
\cmidrule(lr){2-3}\cmidrule(lr){4-5}
Condition & mini & 4o & mini & 4o \\
\midrule
Observation-only                   & 0.27 & 0.25 & 0.24 & 0.17 \\
Point estimate                     & 0.50 & 0.50 & 0.60 & 0.60 \\
Qualitative tags                   & 1.00 & 0.98 & 0.63 & 0.72 \\
Quantitative, means only           & 0.77 & 0.80 & 0.64 & 0.66 \\
$+\ \Delta_{ij}$ only               & 0.77 & 0.98 & 0.63 & 0.70 \\
Uncertainty-aware                  & 0.98 & 0.98 & 0.67 & 0.67 \\
\bottomrule
\end{tabular}

\vspace{.35cm}

\textit{(b) Detection of fingerprint-defined ambiguity}

\vspace{.2cm}
\begin{tabular}{@{}llcc@{}}
\toprule
System (model) & Condition & Ambiguous & Other events \\
\midrule
Idealized (mini) & Qualitative tags & 15/16 & 0/44 \\
Idealized (mini) & $+\ \Delta_{ij}$ only & 1/16 & 1/44 \\
Idealized (mini) & Uncertainty-aware & 16/16 & 1/44 \\
Idealized (4o) & Qualitative tags & 15/16 & 0/44 \\
Idealized (4o) & $+\ \Delta_{ij}$ only & 14/16 & 0/44 \\
Idealized (4o) & Uncertainty-aware & 16/16 & 0/44 \\
Coupled ENSO (mini) & Qualitative tags & 2/8 & 4/82 \\
Coupled ENSO (mini) & $+\ \Delta_{ij}$ only & 0/8 & 0/82 \\
Coupled ENSO (mini) & Uncertainty-aware & 7/8 & 2/82 \\
Coupled ENSO (4o) & Qualitative tags & 7/8 & 10/82 \\
Coupled ENSO (4o) & $+\ \Delta_{ij}$ only & 7/8 & 4/82 \\
Coupled ENSO (4o) & Uncertainty-aware & 7/8 & 1/82 \\
\bottomrule
\end{tabular}

\caption{
Results from the controlled experiments. 
The idealized events are the $60$ with no gate, $Y_1$, $Y_2$, or both gated; the fingerprint rule labels $16$ as ambiguous (the $15$ both-gated events and one $Y_2$-gated event), and $59$ of $60$ labels match the reference. The coupled ENSO experiment contains 90 independent events. GPT-4o-mini and GPT-4o are abbreviated as mini and 4o. (a) Fraction of events for which the LLM attribution agrees with the prescribed intervention label. Six of the eight evidence conditions are shown; the complete results appear in Appendix~\ref{app:ladder-full}.
(b) Number of events diagnosed by the LLM as having indistinguishable mechanisms. The ambiguous group comprises events identified as ambiguous by the fingerprint rule (16 idealized and 8 coupled ENSO events). The remaining 44 and 82 events form the non-ambiguous group. Entries in the latter column are false ambiguity flags relative to the fingerprint rule.
}
\label{tab:ladder-main}
\end{table}

\subsection{An idealized stochastic benchmark}

In the idealized experiment, the mechanism indicators have relatively distinct statistical signatures. When given only the trend and peak of the observed index, GPT-4o-mini and GPT-4o achieve attribution agreement of 0.27 and 0.25, respectively (Table~\ref{tab:ladder-main}a). Providing the deterministic mechanism label raises agreement to 0.50 for both models. These comparisons show the limitations of the restricted observational summary supplied in the first condition, rather than establishing that the full observation record contains no information about the mechanisms.

Qualitative mechanism scores produce much higher agreement, reaching 1.00 for GPT-4o-mini and 0.98 for GPT-4o. Numerical means alone give lower values of 0.77 and 0.80. This difference should not be interpreted as evidence that qualitative information is inherently easier for LLMs to use, since the qualitative and quantitative prompts also differ in their accompanying information (Appendix~\ref{app:ladder-full}).

The effect of the separation statistic differs markedly between the two models. Adding $\Delta_{ij}$ to the numerical means raises GPT-4o agreement from 0.80 to 0.98, with 11 events corrected and none becoming incorrect (exact McNemar test, $p=0.001$). GPT-4o-mini remains at 0.77. The same contrast appears in the ambiguity diagnosis (Table~\ref{tab:ladder-main}b): with numerical separation supplied, GPT-4o identifies 14 of the 16 fingerprint-ambiguous events, compared with only one for GPT-4o-mini.

The uncertainty-aware condition, which supplies the full fingerprint and a verbal interpretation, identifies all 16 ambiguous events under both models. Qualitative tags already identify 15 of 16, however. Thus, the idealized experiment shows that GPT-4o can use numerical separation under these prompts, whereas GPT-4o-mini performs better when the full fingerprint and its verbal conclusion are supplied together. The comparison does not isolate the effect of the verbal conclusion alone. It also shows that a complete numerical description is not always necessary when the driver signatures are relatively easy to distinguish.

\subsection{Mechanism attribution in a coupled ENSO model}

The coupled ENSO experiments illustrate a difficulty of physical attribution: thermocline feedback, zonal advection, and wind bursts interact, so the largest inferred mechanism score need not imply a uniquely resolved explanation. The separation distributions of single- and mixed-intervention events overlap (Appendix~\ref{app:thresholds}). Because the benchmarks also differ in event construction and thresholds, their performance differences cannot be attributed solely to dynamical complexity.

With the observation-only input, attribution agreement is 0.24 for GPT-4o-mini and 0.17 for GPT-4o (Table~\ref{tab:ladder-main}a). The point-estimate condition gives 0.60 for both. Numerical means yield 0.64 and 0.66, while adding $\Delta_{ij}$ changes these values to 0.63 and 0.70. For GPT-4o, the latter comparison corrects six events but makes two previously correct events incorrect; the paired difference is not statistically significant ($p=0.29$). The uncertainty-aware condition gives an agreement of 0.67 for both models. For GPT-4o, qualitative tags yield 0.72, higher than under the uncertainty-aware condition.

Adding score spread or ensemble agreement to numerical means yields no statistically detectable improvement in either system (Appendix~\ref{app:ladder-full}). Neither directly provides the covariance-normalized difference between mechanism scores, but their limited effect under these prompts does not imply that they are generally uninformative.

The coupled ENSO examples illustrate why the distinction between a leading score and a resolved attribution matters. In the mixed thermocline-advection event of Figure~\ref{fig:saicm-diagnostics-app}a (row 4), both mechanisms exceed the support threshold, with mean scores of 2.02 and 1.64. Their separation is only $|\Delta|=0.16$, below the threshold of 1.0. A point estimate favors thermocline, whereas the fingerprint rule assigns a mixed label. The advective score also has substantial ensemble spread ($1.64\pm2.39$); here, \emph{mixed} reflects unresolved ranking under the rule, not proof of equally important physical contributions. In contrast, the thermocline-reference event in Figure~\ref{fig:saicm-diagnostics-app}a (row 1) has two elevated scores, 5.00 for thermocline and 2.28 for wind bursts, but their separation is $|\Delta|=2.79$. The fingerprint therefore favors thermocline alone. A second mixed example in Figure~\ref{fig:uq_blocka} has a larger advective than thermocline mean (2.82 versus 2.46), yet its small separation ($|\Delta|=0.53$) again leads to a mixed fingerprint label. These examples show why neither the highest mean score nor the number of elevated scores is sufficient to determine the fingerprint classification.

The difference between evidence conditions becomes clearer when we examine ambiguity diagnosis rather than overall attribution agreement. Of the 90 coupled ENSO events, eight are ambiguous under the fingerprint rule. With GPT-4o-mini, qualitative tags identify two of these eight events, and the condition adding $\Delta_{ij}$ identifies none (Table~\ref{tab:ladder-main}b). The uncertainty-aware condition identifies seven of eight, while incorrectly flagging two of the remaining 82 events.

GPT-4o identifies seven of eight ambiguous events under each of the three conditions in Table~\ref{tab:ladder-main}b, but the false ambiguity flags differ substantially. Qualitative tags flag 10 of the 82 non-ambiguous events, compared with four when numerical separation is supplied and only one under the uncertainty-aware condition. Relative to qualitative tags, the uncertainty-aware condition eliminates nine false flags without reducing the number of ambiguous events detected (exact McNemar test, $p=0.0039$). The two conditions do not identify exactly the same ambiguous events, however, agreeing on six of the eight cases.

The results distinguish agreement with intervention labels from fidelity to the fingerprint's uncertainty diagnosis. In the coupled experiment, supplying the full fingerprint and its interpretation primarily reduces false ambiguity reports rather than consistently improving intervention-label agreement. Because this condition adds both statistics and a verbal conclusion, their separate effects cannot be identified.

The reliability of these diagnoses also depends on the fingerprint construction. Sensitivity tests revealed that some initially proposed indicators responded strongly to perturbations of other mechanisms, motivating a revised indicator construction (Appendix~\ref{app:indicators}). The ensemble spread also depends on the observation frequency used in data assimilation (Appendix~\ref{app:density}). These checks concern the physical and statistical evidence supplied to the LLM, independently of how well the LLM interprets it.

\subsection{Application to historical ENSO events}

We next examine eight historical El Ni\~{n}o and La Ni\~{n}a events using GODAS SST (Figure~\ref{fig:godas_hov}). Unlike the controlled experiments, these cases have no prescribed perturbation labels. We therefore evaluate each LLM response against the mechanism label derived from the assimilated fingerprint, rather than against an independently established physical attribution.

The uncertainty-aware condition reproduces the fingerprint label for seven of the eight historical events (Table~\ref{tab:godas_events}), with the only mismatch occurring for the 1994 central-Pacific El Niño.

The eight historical cases provide a further check of how an LLM reports model-inferred fingerprints. It does not independently validate the mechanisms assigned to these events. In particular, the high agreement under the uncertainty-aware condition reflects fidelity to the supplied diagnostic conclusion, rather than evidence that the inferred physical attribution is correct.


\section{Discussion and Conclusions}
\label{sec:discussion}

Mechanism attribution in Earth system science involves two distinct questions: what uncertain physical evidence can distinguish, and how faithfully a language model communicates that judgment. In our tests, supplying numerical uncertainty statistics alone did not ensure that the models made this distinction correctly. Pairwise separation improved GPT-4o's intervention-label agreement in the idealized system, but not clearly in the coupled ENSO model. In the latter, both models identified seven of eight fingerprint-ambiguous events when given the full fingerprint and its verbal interpretation. For GPT-4o, false ambiguity reports decreased from 10/82 with qualitative tags to 1/82 with the full condition. Its overall agreement with the imposed intervention labels, however, was slightly lower (0.67 versus 0.72).

These findings suggest separating quantitative mechanism assessment from its communication. Physical and statistical methods can first assess mechanism support and joint uncertainty, while a language model communicates the resulting diagnosis. This does not mean that language models cannot interpret numerical evidence: GPT-4o used the separation statistic effectively in the idealized system. The coupled experiments nevertheless show why it is important to evaluate how accurately a language model reports the uncertainty assessment, rather than relying only on agreement with the imposed intervention labels.

An unresolved fingerprint requires careful interpretation. A small separation indicates that the available fingerprint does not clearly distinguish the mean-supported mechanisms; it does not establish that their causal contributions are equal or that further observations could not resolve the difference. A mixed fingerprint can also arise when one supported mean score has considerable ensemble spread. Similarly, an event with no supported mechanism does not necessarily imply inadequate observations. We did not independently evaluate the diagnosis of insufficient observations through changes to the observing network. Determining which additional measurements would improve attribution therefore remains an open question.

The results also depend on the mechanism indicators and the experimental design. Sensitivity tests show that a mechanism indicator can respond strongly to changes in another physical process (Appendix~\ref{app:indicators}), making validation of the indicators essential. The ensemble covariance entering $\Delta_{ij}$ also depends on the assimilation configuration, and posterior calibration against the known reference states was not established (Appendix~\ref{app:density}). The intervention-based reference labels identify the deliberately amplified processes but do not provide a unique causal decomposition of each stochastic event. Agreement with these labels combines the limitations of the inferred fingerprint and the LLM interpretation, which are not separately resolved by the end-to-end scores. For the historical ENSO cases, agreement with the model-inferred fingerprints measures communication fidelity rather than independently validated physical attribution. The historical indicators can also reflect model-tendency residuals and use a different normalization from the controlled experiments. Furthermore, the quantitative conditions share information beyond the individual statistics being compared, and the uncertainty-aware condition combines additional statistics with a verbal interpretation. The separate contribution of the verbal statement therefore remains unresolved. Our evaluation is also limited to two related language models and one response per event and condition. More recent or differently trained models may apply the numerical decision rule more reliably, so these findings should not be assumed to hold across LLMs. Regardless of model capability, validating the physical uncertainty diagnostics remains distinct from testing how accurately a language model communicates them.

Although ENSO is the physical testbed here, assessing uncertainty and communicating its implications are separate tasks whenever language models explain uncertain Earth system diagnostics. Probabilistic fingerprints require mechanism-related diagnostics and their joint uncertainty, not a particular model or assimilation method. Whether the response patterns persist in other geophysical settings and newer language models remains to be tested.

\section*{Open Research}
	
Code, event libraries, and full prompts and responses, including model identifiers, are available in \href{https://github.com/pbehnoud/Uncertainity_Aware_LLM.git}{GitHub}. GODAS data were obtained from NOAA \citep{behringer2004}. 
	
	\section*{Acknowledgments}
The research of N.C. is funded by the Office of Naval Research N00014-24-1-2244 and the Army Research Office W911NF-23-1-0118. P.B. is partially supported as a research associate under the second grant.
	
	\bibliographystyle{plainnat}
	\bibliography{references}
	
	\clearpage
	\appendix
	\renewcommand{\thesection}{S\arabic{section}}
	\renewcommand{\thefigure}{S\arabic{figure}}
	\renewcommand{\thetable}{S\arabic{table}}
	\setcounter{figure}{0}
	\setcounter{table}{0}
\section*{Supporting Information}
	\addcontentsline{toc}{section}{Supporting Information}

\section{Conceptual model: indicators and separation}
\label{app:concept}

This section shows how the mechanism indicators and the separation statistic label events in the conceptual model. The posterior ensemble is taken as given, and the question is whether the fingerprint computed from it identifies the gated driver, reports two drivers when both are gated, and reports none when no driver is gated.
Figure~\ref{fig:c1_events} shows the posterior and the fingerprint of one event of each type, Figure~\ref{fig:c1_scores} the scores of all events, and Figure~\ref{fig:c1_sep} how the separation decides the label in the cases closest to the thresholds.

\subsection{Model, events, and fingerprint}
\label{app:c1_model}

Three independent Ornstein--Uhlenbeck drivers $Y_1,Y_2,Y_3$ force an observed variable $X$,
\begin{align}
dY_k &= -\lambda_k Y_k\,dt + \sigma_k(t)\,dW_k, \qquad k=1,2,3,\\
dX &= \Big(\textstyle\sum_k Y_k - \kappa X\Big)dt + \sigma_X\,dW_X,
\end{align}
with $\lambda=(0.1,\,0.5,\,0.25)$, $\kappa=0.3$ and $\sigma_X=0.1$.
The drivers differ in their memory: $Y_1$ is the slowest (correlation time $1/\lambda_1=10$), $Y_2$ is the fastest ($2$), and $Y_3$ lies in between ($4$). These are generic stochastic drivers, not representations of particular geophysical mechanisms. Outside events, $\sigma_k=\sigma_0=0.3$. In an event window, the noise amplitude of a gated driver is raised to $\sigma_k(t)=\sigma_0+1.6\,a(t)$, where $a(t)$ is a Gaussian pulse that equals $1$ at the center of the window and $0.32$ at its edges. Only $X$ is observed every $0.8$ time units, with error standard deviation $0.3$. 

The indicator of driver $k$ is $S_k=\int|Y_k|\,dt$ over the attribution window, computed for each member, and standardized as $Z_k=(S_k-\mu_k)/\sigma_{{\rm clim},k}$. Here $\mu_k$ and $\sigma_{{\rm clim},k}$ are the mean and standard deviation of $S_k$ over $500$ windows of the same length drawn from a $4000$-unit run without gates. $Z_k$ therefore measures how far the activity of driver $k$ in the window lies above its ordinary range. The fingerprint is $\bar Z_k$, $\sigma_{Z,k}$ and the separation $\Delta_{ij}$ (Section~\ref{sec:framework}), and the labels use $z_{\rm thr}=1$ and $\delta_{\rm thr}=2.5$.

Each event is an independent run: a burn-in without gates followed by the attribution window (Table~\ref{tab:c1_lib}). There are five event types: with no gate, $Y_1$, $Y_2$, both, or $Y_3$ gated. The LLM experiments use the $60$ events without the $Y_3$-gated type.

\begin{table}[htbp]
\centering
\small
\caption{Events of the conceptual model. The label is correct when it is none for events without a gate, single with the gated driver for single-driver events, and mixed with the two gated drivers for events with both. Score-separation column counts the number of correct labels using the indicator label and separation $\Delta$. The last column is for the case that names the driver with the largest $\bar Z$.}
\label{tab:c1_lib}
\begin{tabular}{@{}lccccc@{}}
\toprule
Gated driver(s) & burn-in & window & events & Score-separation correct & Score correct \\
\midrule
none & $20$ & $20$ & $15$ & $15$ & $0$ \\
$Y_1$ & $30$ & $25$ & $15$ & $15$ & $15$ \\
$Y_2$ & $30$ & $25$ & $15$ & $14$ & $15$ \\
$Y_3$ & $30$ & $25$ & $15$ & $15$ & $15$ \\
$Y_1$ and $Y_2$ & $30$ & $40$ & $15$ & $15$ & $0$ \\
\midrule
total & & & $75$ & $74$ & $45$ \\
\bottomrule
\end{tabular}
\end{table}

\subsection{One event of each type}
\label{app:c1_events}

\begin{figure}[htbp]
\centering
\includegraphics[width=\textwidth]{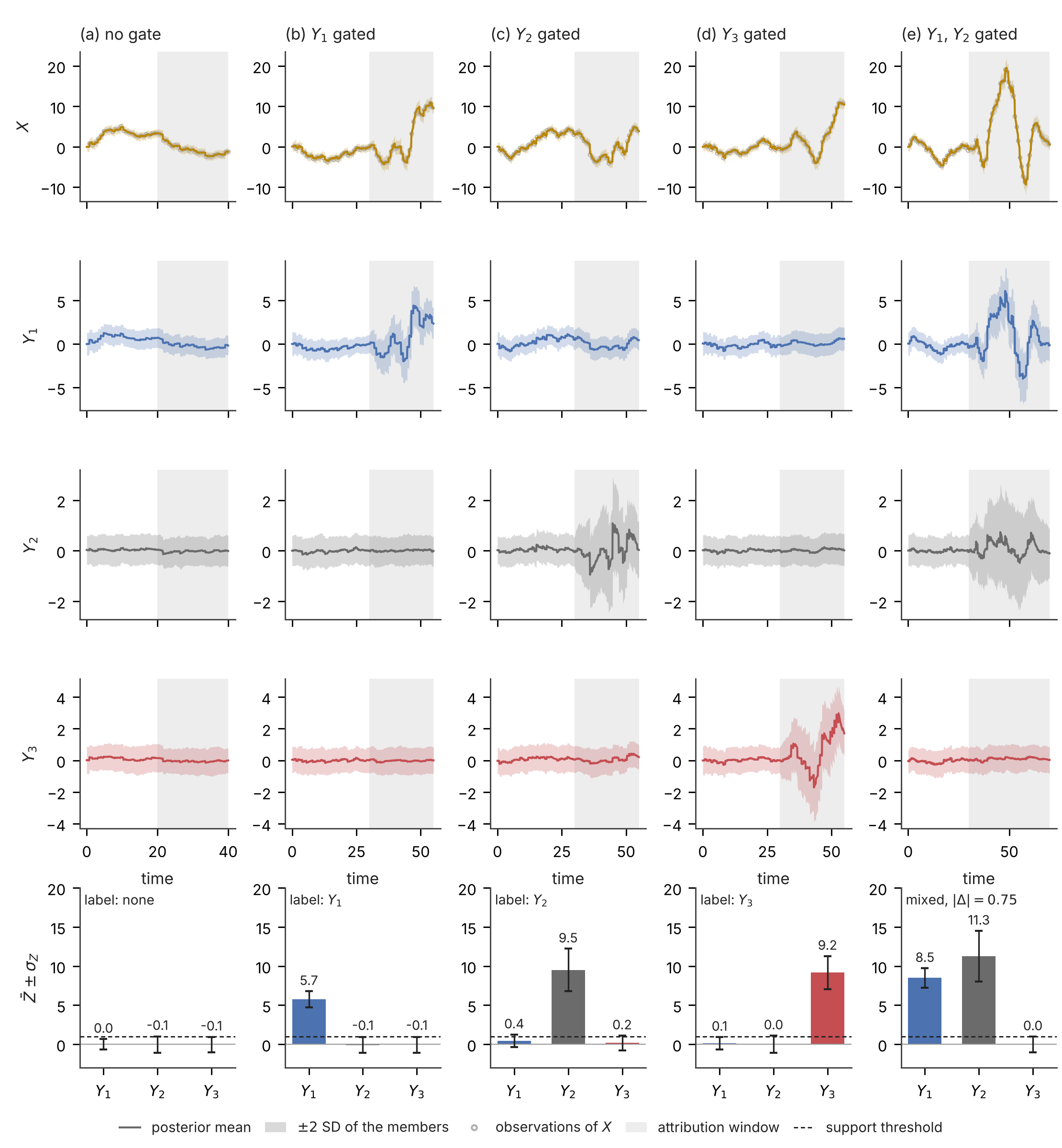}
\caption{Posterior and fingerprint of one event of each type: (a) no gate, (b) $Y_1$ gated, (c) $Y_2$ gated, (d) $Y_3$ gated, (e) $Y_1$ and $Y_2$ gated. Rows 1--4: posterior mean (line) and $\pm2$ standard deviations of the ensemble members (band) of the observed variable $X$, with its observations as open circles, and of the three drivers. Grey shading marks the attribution window. Each row has one vertical scale. Row 5: fingerprint $\bar Z\pm\sigma_Z$ of the $Y_1$, $Y_2$, and $Y_3$ indicators, with the support threshold (dashed) and the resulting label.}
\label{fig:c1_events}
\end{figure}

The indicators respond to what happens to each driver inside the window (Figure~\ref{fig:c1_events}).

\paragraph{No gate (a).}
All three drivers keep the same behavior before and inside the window, and the integral of $|Y_k|$ stays within its ordinary range. The scores are $-0.01$, $-0.06$ and $-0.06$, none is supported, and the label is none.

\paragraph{One driver gated (b--d).}
In each case only the gated driver changes in the window: its posterior shows larger excursions and a wider band, and its integrated activity rises well above the climatological range. The corresponding score is $5.74\pm1.03$ for $Y_1$ (b), $9.50\pm2.72$ for $Y_2$ (c) and $9.15\pm2.13$ for $Y_3$ (d), while the two other scores stay between $-0.10$ and $0.44$. Only the gated driver is supported, so each event is labeled single with the right driver. The spread differs between the drivers. The posterior band of $Y_1$ widens only moderately, whereas that of $Y_2$ and $Y_3$ widens strongly in the window. The members then disagree more on how active the driver was, and $\sigma_Z$ is larger relative to $\bar Z$: the median of $\sigma_Z/\bar Z$ is $0.21$ for the $Y_1$ events and $0.28$ for both the $Y_2$ and $Y_3$ events.

\paragraph{Two drivers gated (e).}
Both $Y_1$ and $Y_2$ change in the window. Both scores are supported ($8.48\pm1.24$ and $11.26\pm3.24$), and the $Y_3$ score stays at $-0.03$. The $Y_2$ score is the larger, but its difference from the $Y_1$ score, $2.78$, is only $0.75$ standard deviations of the difference across members ($\Delta=-0.75$). The label is mixed: both drivers are active, and the posterior does not establish that one dominates.

\subsection{Each indicator responds to its own driver}
\label{app:c1_specific}

\begin{figure}[htbp]
\centering
\includegraphics[width=\textwidth]{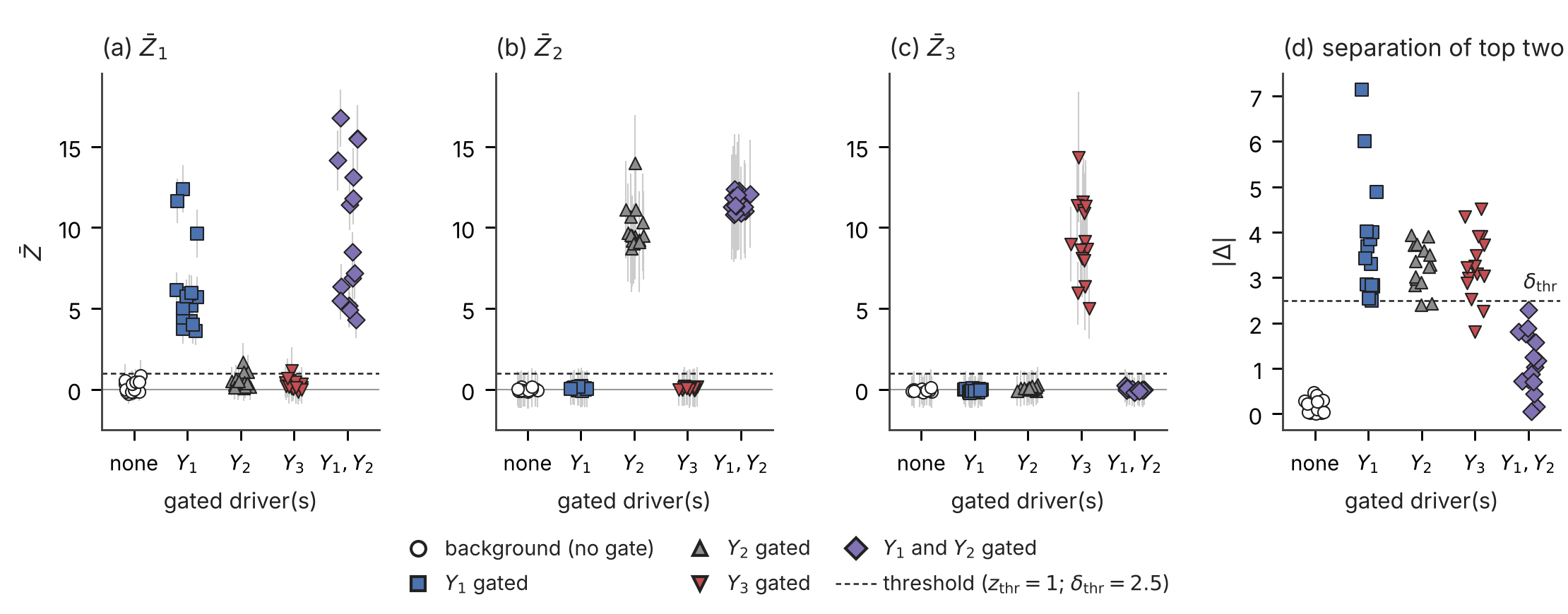}
\caption{Fingerprints of all $75$ events. (a--c) Mean score $\bar Z$ of the $Y_1$, $Y_2$ and $Y_3$ indicators, grouped by the gated driver(s), with $\pm\sigma_Z$ as grey bars; the dashed line is $z_{\rm thr}=1$. (d) Separation $|\Delta|$ between the two highest scores; the dashed line is $\delta_{\rm thr}=2.5$.}
\label{fig:c1_scores}
\end{figure}

\begin{table}[htbp]
\centering
\small
\caption{Mean score $\bar Z$ of each indicator over the $15$ events of each type, with the range over the events in brackets.}
\label{tab:c1_spec}
\begin{tabular}{@{}lccc@{}}
\toprule
Gated driver(s) & $Y_1$ & $Y_2$ & $Y_3$ \\
\midrule
none & $0.15$ $[-0.29,\,0.85]$ & $-0.02$ $[-0.15,\,0.14]$ & $-0.05$ $[-0.19,\,0.11]$ \\
$Y_1$ & $\mathbf{6.22}$ $[3.60,\,12.38]$ & $0.05$ $[-0.10,\,0.19]$ & $-0.06$ $[-0.22,\,0.07]$ \\
$Y_2$ & $0.53$ $[0.09,\,1.66]$ & $\mathbf{10.00}$ $[8.72,\,13.97]$ & $0.06$ $[-0.10,\,0.28]$ \\
$Y_3$ & $0.29$ $[-0.16,\,1.15]$ & $0.05$ $[-0.06,\,0.19]$ & $\mathbf{9.28}$ $[4.99,\,14.34]$ \\
$Y_1$ and $Y_2$ & $\mathbf{9.80}$ $[4.32,\,16.75]$ & $\mathbf{11.50}$ $[10.81,\,12.36]$ & $-0.01$ $[-0.18,\,0.24]$ \\
\bottomrule
\end{tabular}
\end{table}

Across all events, each indicator responds most strongly when its own driver is gated (Figure~\ref{fig:c1_scores}a--c, Table~\ref{tab:c1_spec}). The scores of gated drivers range from $3.6$ to $16.8$, and their distance above the threshold, $(\bar Z-1)/\sigma_Z$, has a median of $3.4$ (range $2.0$ to $9.0$). No gated driver comes close to losing its support. Without a gate, all scores lie between $-0.29$ and $0.85$, so the threshold $z_{\rm thr}=1$ produces no false support in the $15$ background events. Scores of drivers that are not gated stay near zero in the other types as well, with one exception: the $Y_1$ score rises slightly when a faster driver is gated. It exceeds the threshold in $3$ of the $15$ $Y_2$ events (up to $1.66$) and in $1$ of the $15$ $Y_3$ events ($1.15$). Part of the extra activity in the window is taken up by the slow driver. In these four events, the gated driver still has by far the larger score, and the separation decides whether the second score matters (Section~\ref{app:c1_sep}).

\subsection{The separation decides the label}
\label{app:c1_sep}

\begin{figure}[htbp]
\centering
\includegraphics[width=\textwidth]{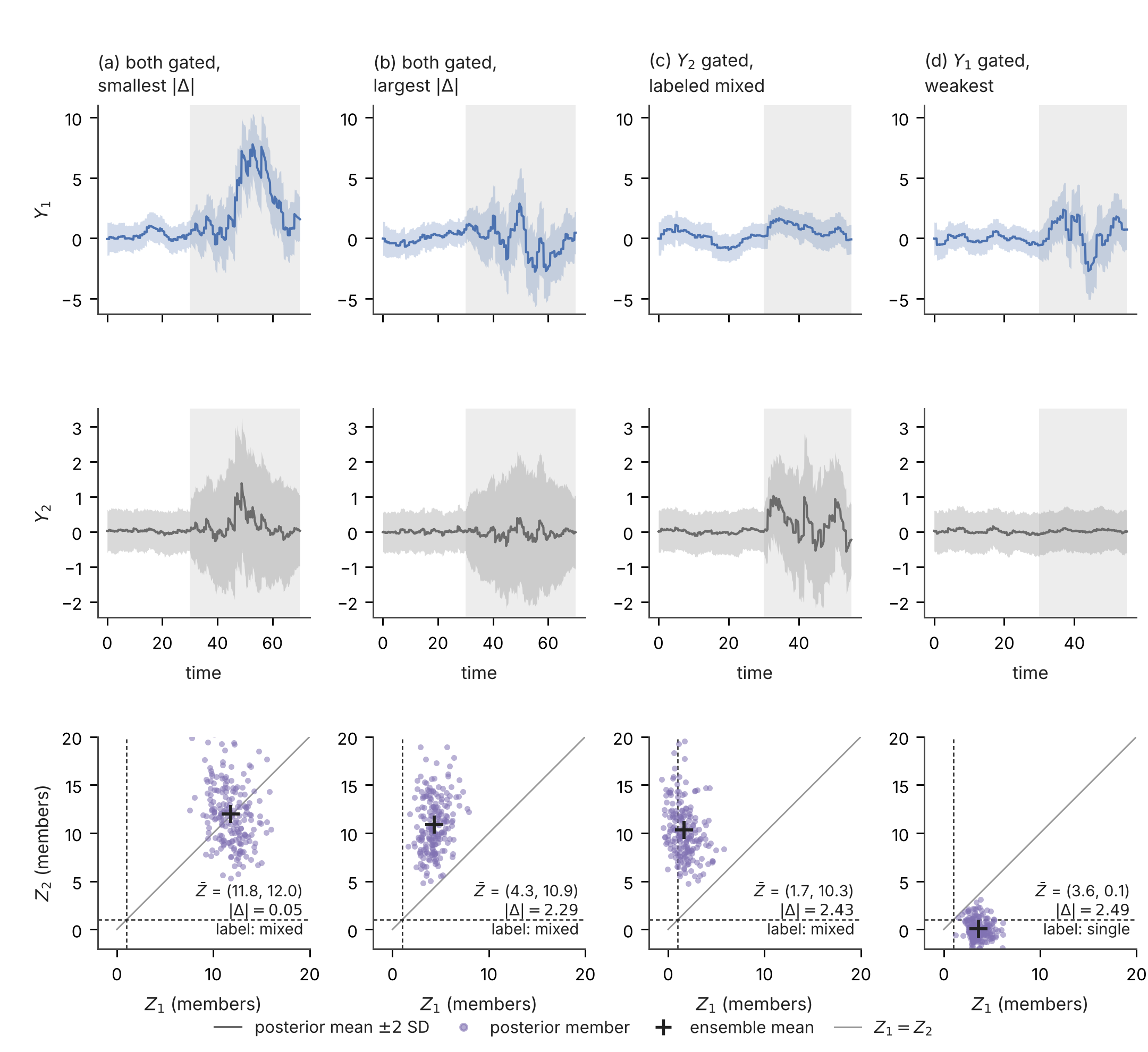}
\caption{Four events close to the thresholds. (a) $Y_1$ and $Y_2$ gated, smallest $|\Delta|$. (b) $Y_1$ and $Y_2$ gated, largest $|\Delta|$. (c) $Y_2$ gated, labeled mixed. (d) $Y_1$ gated, weakest score. Top rows: posterior mean and $\pm2$ standard deviations of $Y_1$ and $Y_2$. Bottom row: $Y_1$ and $Y_2$ scores of the ensemble members (dots) and their mean (cross); the dashed lines mark the support threshold, and the grey line denotes $Z_1=Z_2$ (labeled $Z_{\rm Th}=Z_{\rm Wi}$ in the existing figure).}
\label{fig:c1_sep}
\end{figure}

The lower row of Figure~\ref{fig:c1_sep} shows the scores of the individual members, which is what the fingerprint summarizes. The position of the cloud relative to the thresholds gives the support, and its position relative to the diagonal gives the separation. In the events with two gated drivers, the member scores are nearly uncorrelated (mean correlation $-0.15$), so $\Delta$ is close to the difference of the two mean scores divided by $\sqrt{\sigma_{Z,1}^2+\sigma_{Z,2}^2}$.

\paragraph{Two gated drivers.}
In Figure~\ref{fig:c1_sep}a the scores are $11.8\pm1.6$ and $12.0\pm3.3$. The cloud lies across the diagonal, $|\Delta|=0.05$, and nothing in the posterior favors one driver. In event~$4$ (b) the scores are $4.3\pm1.1$ and $10.9\pm2.9$. The cloud lies on one side of the diagonal, but the difference of $6.6$ is only $2.29$ standard deviations of the difference, below $\delta_{\rm thr}$, and the event is labeled mixed.
These two events bound the both-gated set: $|\Delta|$ ranges from $0.05$ to $2.29$, and $7$ of the $15$ events are below $1$ (Figure~\ref{fig:c1_scores}d).
In all $15$ events both scores are supported and $|\Delta|<\delta_{\rm thr}$, so all are labeled mixed. A point estimate would name a single driver in each: driver $Y_2$ in $10$ events and driver $Y_1$ in $5$.

\paragraph{One gated driver.}
For single-driver events, the separation of the top two scores ranges from $2.49$ to $7.14$ ($Y_1$), $2.41$ to $3.93$ ($Y_2$), and $1.82$ to $4.52$ ($Y_3$).
Five of these $45$ events have a separation below $\delta_{\rm thr}$. Four are still labeled single because their second score is not supported, so the separation does not apply. Figure~\ref{fig:c1_sep}d is one of them: its $Y_1$ score of $3.6\pm0.9$ is the weakest of its type, and $|\Delta|=2.49$, but the second score is $0.06$.
The fifth is the only mislabeled event of the $75$ (Figure~\ref{fig:c1_sep}c). Here, the $Y_1$ score of a $Y_2$-gated event is $1.66\pm1.22$, just above the threshold, and its separation from the $Y_2$ score ($10.31\pm3.00$) is $2.43$. The cloud lies well off the diagonal but extends below $Z_1=1$, and the rule reports both drivers.
The other three single-driver events with a supported second score have separations of $2.90$ and $3.91$ ($Y_2$) and $2.52$ ($Y_3$, just above the threshold), and are labeled single.

\paragraph{Without the separation.}
A point estimate, which names the driver with the largest $\bar Z$, is correct for all $45$ single-driver events but cannot return mixed or none, so it fails on all $30$ events with no gate or two gates ($45$ of $75$, Table~\ref{tab:c1_lib}).
The support threshold handles the events without a gate. The separation is what distinguishes the $15$ events with two active drivers from the $45$ with one: between them, the two quantities give $74$ correct labels.
\section{Experimental systems and controlled perturbations}
\label{app:systems}

\textbf{The conceptual model} consists of three independent Ornstein--Uhlenbeck drivers with distinct relaxation timescales that force a common observed variable. Within a prescribed event window, the forcing amplitude of one or more selected drivers is increased. The unmodified components remain active; a lack of imposed perturbation should not be interpreted as absence of physical dynamics.
	
\begin{figure}[htbp]
		\centering
		\includegraphics[width=1\linewidth]{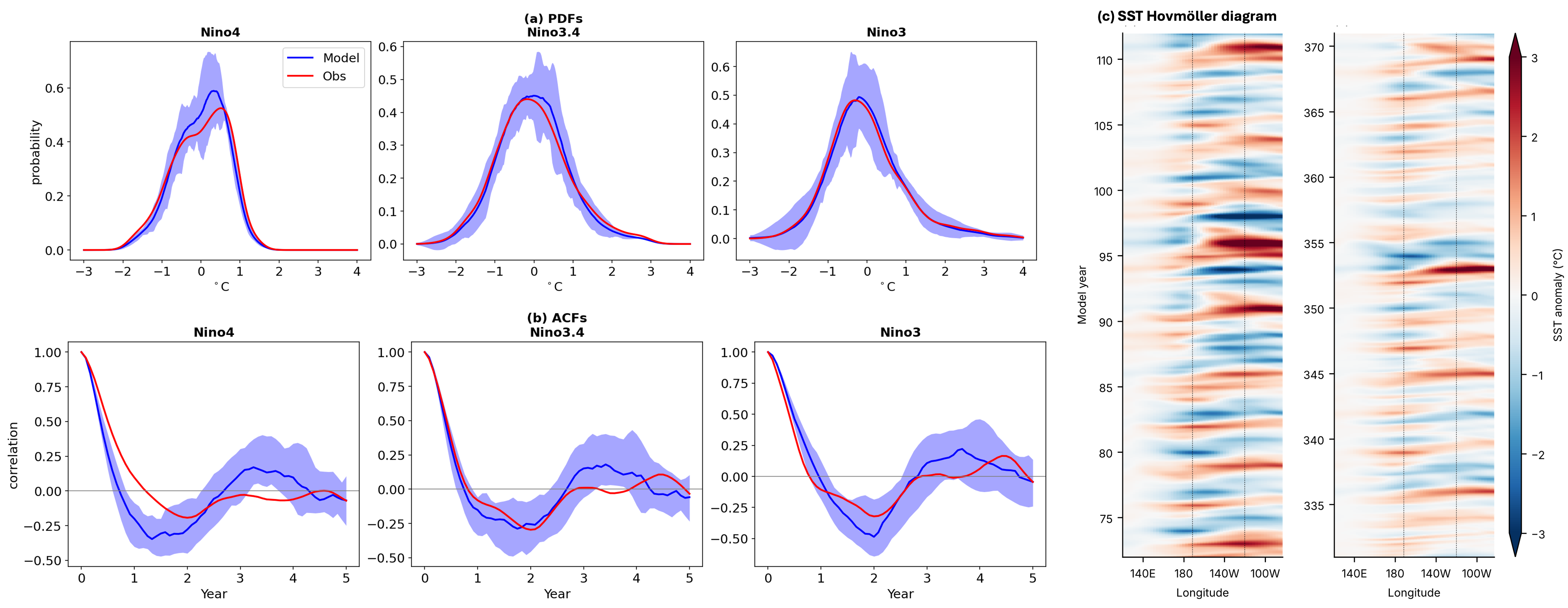}
		\caption{Comparison of the PDFs (a) and autocorrelation functions (b) of SST anomalies between the ICM and GODAS in the Ni\~{n}o regions. The blue shaded area spans two standard deviations of the statistics across model segments of the same length as the observational record.(c)  the sea surface temperature Hovm\"oller  over 40-year intervals showing the ICM's capacity for capturing ENSO variablity.}
		\label{fig:stat}
	\end{figure}

\textbf{The ICM} resolves SST and the leading equatorial Kelvin and Rossby mode amplitudes along a zonal grid spanning the equatorial Pacific ($120^\circ$--$280^\circ$E) \citep{chen2023}, and supports a spatial diversity diagnostic (central- versus eastern-Pacific; \citealp{kug2009,renjin2011}) computed from peak-time Ni\~{n}o3 and Ni\~{n}o4 indices. Figure~\ref{fig:stat} compares its SST statistics with GODAS reanalysis data. Representative ICM space--time fields are already shown in the Hovm\"oller diagrams of Figures~\ref{fig:saicm-diagnostics-app}, \ref{fig:uq_blocka}, and~\ref{fig:icm_hov_1}; these illustrate simulated dynamics but are not an independent validation against historical observations.

\textbf{Event construction.} For each system, an event consisted of an independent simulation with an initial background period without imposed amplification, followed by an event window in which one or two selected processes are strengthened. Reference labels are defined by the imposed perturbations, whereas the response of the observed quantity (SST in the coupled ENSO model) is illustrated by paired perturbed and unperturbed simulations. For this comparison, both simulations start from the same initial conditions and use the same realization of stochastic forcing. 

Reference labels are defined by construction: an event is labeled by the set of processes whose gates were switched on, independently of the simulated response of the observed variable, so no response threshold is applied (the CP/EP/neutral SST class, with a $0.3$ K neutral threshold, is a separate descriptor). Each gate multiplies the corresponding coefficient (thermocline, advection) or amplitude (wind bursts) by $1+A\,a(t)$, where $a(t)$ is a Gaussian pulse centered in the event window with standard deviation one third of the window length. In the ICM, $A=1.5$, $1.5$ and $4.0$ for thermocline, wind-burst and advective gates, with windows of 3, 1.5 and 3 years (4 years for mixed events) after a 15-year burn-in. In the conceptual model, the driver noise is $\sigma_k=\sigma_0+1.6\,\rho_k\,a(t)$ with $\sigma_0=0.3$. Mixed events switch on two gates simultaneously and follow no separate labelling rule, and all events generated per type were retained. The ICM's ENSO statistics were compared with GODAS in an unperturbed control run (Figure~\ref{fig:stat}).

\section{Mechanism indicators and sensitivity tests}
\label{app:indicators}

For each event and each ensemble member, an indicator $S_k$ is computed for every mechanism $k$ over the event window $[t_0,t^*]$. Each score is aggregated over the whole window, and in the ICM, it is also pooled over all grid points. It is then standardized as $Z_k=(S_k-\mu_k)/\sigma_k$, where $\mu_k$ and $\sigma_k$ are the mean and standard deviation of $S_k$ over windows of the same length drawn from an ungated control run (500 windows from a run of $4000$ time units for the conceptual model; 400 windows from a $200$-year run for the ICM). The fingerprint of an event consists of the ensemble mean $\bar Z_k$, the ensemble standard deviation $\sigma_{Z_k}$, and the separation $\Delta_{ij}$ of Section~\ref{sec:framework}. In the conceptual model, $S_k=\int_{t_0}^{t^*}|Y_k|\,dt$, where $Y_k$ is the $k$-th driver. In the ICM, the wind-burst indicator is $S_{\rm wind}=\int_{t_0}^{t^*}|a_p|\,dt$, where $a_p$ is the wind-burst amplitude. The thermocline and advective indicators are the gate amplitudes recovered by the regression described below. They are signed and are standardized against the same statistic computed on the control run.

An indicator for mechanism $k$ should respond primarily to the parameter controlled by that mechanism; however, this expectation can fail without producing any visible error in the output. We identified one such case in each system. In the conceptual model, an indicator measured a driver whose statistical properties were unaffected by the corresponding gate. In the ICM, indicators based on equatorial wave amplitudes responded to the gated SST-tendency coefficients only indirectly through the coupled feedback loop.

	\begin{table}[htbp]
		\centering
		\begin{tabular}{@{}lccc@{}}
			\toprule
			gate $\backslash$ indicator & thermocline & wind burst & advective \\
			\midrule
			thermocline & $+40.0\%$ & $-22.9\%$ & $+5.5\%$ \\
			wind burst & $+74.9\%$ & $+89.9\%$ & $+65.6\%$ \\
			advective & $+16.4\%$ & $+11.2\%$ & $+17.6\%$ \\
			\bottomrule
		\end{tabular}
		\caption{Sensitivity of the original ICM indicators (common-random-number fractional change in $\bar Z_k$ per gate). The advective gate moves its own indicator only slightly more than the thermocline indicator ($+17.6\%$ against $+16.4\%$), and the wind-burst gate moves all three indicators; row correlation $+0.74$.}

		\label{tab:sens-app}
	\end{table}
	
The thermocline and advective gates amplify coefficients of the SST-tendency equation, whereas the amplitude-based indicators measure the wave fields that respond to them only through the coupled feedback. Because SST carries no additive noise in the ICM, its tendency is determined exactly by the model state. The regression therefore differentiates the SST tendency, subtracts the tendency predicted with the baseline coefficients, and projects the residual onto the sum and difference of the basis fields associated with the two gated coefficients. The projection is a least-squares fit for each member over time and space, and its coefficients are the fractional gate amplitudes. The recovered amplitudes were $+1.13/-0.11$ for the thermocline-gated case, $+0.07/+2.19$ for the advective-gated case, and $+0.99/+1.75$ when both mechanisms were gated, compared with $+0.03/+0.03$ in the ungated case. The row correlation of Table~\ref{tab:sens-app} consequently changed from $+0.74$ to $-0.62$.

The wind-burst gate raised all three original indicators (Table~\ref{tab:sens-app}) because wind bursts drive the Kelvin and Rossby amplitudes directly. The wind-burst indicator is retained in its amplitude-based form because it measures the stochastic wind-burst driver itself, and it responds most strongly to its own gate ($+89.9\%$, against $-22.9\%$ and $+11.2\%$ for the thermocline and advective gates).

With the regression-based indicators, the accuracy of the deterministic classifier (the mechanism with the largest $\bar Z_k$) on the ICM, measured against the reference labels, increased from $0.42$ to $0.55$, and its accuracy for \texttt{advective\_only} increased from $0$ to $0.62$. These values were obtained during indicator development. The final experiments use a monthly observation interval, which changes the accuracy (Appendix~\ref{app:density}), so the accuracies reported in the main text differ from these values and are those of Table~\ref{tab:ladder-main}.

\section{Fingerprint decision rules and threshold calibration}
	\label{app:thresholds}

The separation statistic uses the covariance of paired scores from the same ensemble member,
\[
\Delta_{ij}=\frac{\bar Z_i-\bar Z_j}{\sqrt{\sigma_{Z_i}^2+\sigma_{Z_j}^2-2\,\mathrm{Cov}(Z_i,Z_j)}}.
\]
Its denominator is the standard deviation of the memberwise score difference, not the standard error of the estimated mean difference. Omitting positive covariance tends to underestimate separation; omitting negative covariance tends to overestimate it. Either change may alter a decision when an event lies near a threshold.

An indicator is supported when $\bar Z_k>z_{\rm thr}$. When no score is supported, the fingerprint assigns \emph{none}; when one is supported, it assigns that mechanism. For multiple supported scores, the rule compares the two largest mean scores and assigns a \emph{mixed} fingerprint label if their absolute separation falls below $\delta_{\rm thr}$. This is an operational category, not an independently measured causal mixture.

The support threshold $z_{\rm thr}$ is set from background events, where each exceedance is a false support flag relative to the background reference. This pooled criterion does not automatically control the probability of at least one false flag across all candidate mechanisms. We use $z_{\rm thr}=1$ in both systems. On the ICM it is calibrated on the ICM event library as the $90\%$ quantile of the pooled mean scores on the background events, which fixes the false support rate at $10\%$ and uses no information from the perturbed events. On the conceptual model it is one standard deviation of the standardized background scores: in the $15$ background events no mean score exceeds $0.85$, so no indicator is supported, and about $15\%$ of the individual member scores exceed $1$.

The separation $\Delta_{ij}$ matters only when two indicators are supported. In the $75$ conceptual-model events, the four single-driver events with two supported indicators have $|\Delta_{ij}|\in[2.43,3.91]$, and the $15$ mixed events, all with two supported indicators, have $|\Delta_{ij}|\in[0.05,2.29]$. The threshold $\delta_{\rm thr}=2.5$ lies just above the smallest single-driver value, so $74$ of the $75$ events are labeled as prescribed; the exception is a $Y_2$-gated event with $|\Delta_{ij}|=2.43$ (Appendix~\ref{app:concept}). On the ICM, $\delta_{\rm thr}=1.0$ is a standardized difference of one under the defined separation measure. These thresholds define the fingerprint labels against which the LLM conditions are evaluated.

\section{Ensemble assimilation and observation-frequency sensitivity}
\label{app:density}

In the controlled coupled ENSO experiments, a stochastic ensemble Kalman filter with 150 members, perturbed observations and relaxation to the prior spread is used to infer model variables from SST \citep{evensen1994,burgers1998,whitaker2012}. The observation network consists of SST measurements at the model's 56 longitudinal grid points, with an error standard deviation of $r=0.3^\circ$C. Kelvin and Rossby amplitudes and wind forcing are not directly observed in this experiment. The final benchmark adopts monthly assimilation (60 model steps per month). The filter updates the Kelvin and Rossby amplitudes, SST, the wind-burst amplitude, and the decadal index jointly. Observations are perturbed with standard deviation $r$ at each analysis step.

With SST observed at all $56$ grid points every $8$ steps, the ratio of posterior spread to observation-error standard deviation was $0.17$, and $82\%$ of events had ensemble agreement above $80\%$. Increasing observation error to $r=2.0$ gave a ratio of $0.08$ in the reported experiment. Because this ratio uses $r$ as its denominator, a lower ratio need not imply a better- or worse-calibrated ensemble. At fixed $r=0.3^\circ$C, decreasing observation frequency increased the absolute posterior spread as shown below:
	
	\begin{table}[htbp]
		\centering
		\begin{tabular}{@{}rrcc@{}}
			\toprule
			Obs.\ interval (steps) & analysis cycles & posterior spread & ratio to $r$ \\
			\midrule
			8 & 1710 & 0.050 & 0.17 \\
			25 & 547 & 0.080 & 0.27 \\
			60 & 228 & 0.122 & 0.41 \\
			150 & 91 & 0.227 & 0.76 \\
			400 & 34 & 0.414 & 1.38 \\
			\bottomrule
		\end{tabular}
		\caption{Posterior spread against observation interval at $r=0.3^\circ$C, for a single event (mixed thermocline--advective, first realization). The spread is the ensemble standard deviation of SST, averaged over the 56 grid points and the full run (burn-in and event window), in $^\circ$C. The number of analysis cycles is the number of assimilation times, and each cycle assimilates SST at all 56 grid points.}
		\label{tab:density-app}
	\end{table}

One model month is $60$ steps, so an interval of $60$ matches the monthly reporting of SST indices; it gives a ratio of $0.41$ and was accompanied by a reported change in deterministic classifier agreement from $0.55$ to $0.62$. The comparison does not establish filter calibration; this would require spread-error or coverage diagnostics against the known reference states.


\section{LLM evidence conditions and evaluation}
	\label{app:ladder-full}

The two language models are GPT-4o-mini (a smaller and less computationally expensive model) and GPT-4o (a higher-capability model in the same family). In the archived conceptual-model prompts, the descriptive labels thermocline-like, wind-like, and advection-like are mnemonics for $Y_1$, $Y_2$, and $Y_3$, respectively; these generic drivers are not geophysical mechanisms. System instructions, generation settings and structured output format remain fixed across conditions. Only supplied evidence varies; this design does not isolate the effects of all individual prompt fields, since the qualitative and quantitative formats also differ in their contextual information.

The eight conditions are observation-only, point estimate, qualitative tags, numerical means, means plus score spread, means plus pairwise separation, means plus ensemble agreement, and the full fingerprint together with a verbal conclusion. The quantitative conditions retain ensemble-averaged softmax support weights, and the coupled ENSO versions also contain SST diversity diagnostics (ensemble vote agreement and regional breakdown). Thus, the means-only condition already contains some ensemble information, and the comparisons do not isolate the value of uncertainty as a whole.

Each event was submitted once per model and condition as a single chat-completion request (temperature $0.3$, maximum $500$ output tokens) with a fixed system prompt and a fixed structured-output template; event identifiers were replaced by opaque codes. The models were accessed through the OpenAI API under the names \texttt{gpt-4o} and \texttt{gpt-4o-mini}. All fingerprint conditions state that a mechanism is supported only if its mean standardized score exceeds $1.0$, and conditions that include $\Delta_{ij}$ state that $|\Delta_{ij}|<2.5$ (conceptual model) or $<1.0$ (coupled ENSO model) indicates mechanisms that are hard to distinguish, but only if both are supported. In the qualitative condition, $\bar Z_k$ is reported as high ($>2$), moderate ($>1$), weak ($>0$) or negligible. Table~\ref{tab:prompt-fields} lists the fields supplied in each condition. A response was scored from its \texttt{DOMINANT\_MECHANISM} field alone: correct if it equals the reference mechanism, \texttt{mixed} for mixed events, or \texttt{none} for background events. The named pair was not required to match; requiring it changes two entries of Table~\ref{tab:ladder-full-app} by $0.01$. Table~\ref{tab:prompt-fields} provides an overview of the pipeline, and the complete prompts and responses are available at \href{https://github.com/pbehnoud/Uncertainity_Aware_LLM.git}{GitHub}.

\begin{table}[htbp]
\centering\scriptsize
\begin{tabular}{@{}lcccccccc@{}}
\toprule
Field & Obs. & Point & Tags & Means & $+\sigma_Z$ & $+\Delta_{ij}$ & $+$agr. & Full \\
\midrule
Observed index (trend, peak value and time) & \checkmark & & & & & & & \\
Observation context (variable, noise, interval, ensemble size) & & \checkmark & & \checkmark & \checkmark & \checkmark & \checkmark & \checkmark \\
Mechanism chosen by the classifier & & \checkmark & & & & & & \\
Qualitative tag per mechanism & & & \checkmark & & & & & \\
Mean score $\bar Z_k$ and support threshold & & & & \checkmark & \checkmark & \checkmark & \checkmark & \checkmark \\
Relative support weight & & & & \checkmark & \checkmark & \checkmark & \checkmark & \checkmark \\
Posterior standard deviation $\sigma_{Z_k}$ & & & & & \checkmark & & & \checkmark \\
Pairwise separation $\Delta_{ij}$ and threshold & & & & & & \checkmark & & \checkmark \\
Ensemble agreement & & & & & & & \checkmark & \checkmark \\
Verbal summary of the conclusion & & & & & & & & \checkmark \\
\bottomrule
\end{tabular}
\caption{Fields supplied in each evidence condition (conceptual model). In the coupled ENSO prompts, the quantitative conditions also contain the SST diversity classification and its regional breakdown.}
\label{tab:prompt-fields}
\end{table}

	\begin{table}[htbp]
		\centering
		\small
		\begin{tabular}{@{}lcccc@{}}
			\toprule
			& \multicolumn{2}{c}{Idealized ($n=60$)} & \multicolumn{2}{c}{Coupled ENSO ($n=90$)} \\
			\cmidrule(lr){2-3}\cmidrule(lr){4-5}
			Condition & mini & 4o & mini & 4o \\
			\midrule
			Observation-only          & 0.27 & 0.25 & 0.24 & 0.17 \\
			Point estimate            & 0.50 & 0.50 & 0.60 & 0.60 \\
			Qualitative tags          & 1.00 & 0.98 & 0.63 & 0.72 \\
			Quantitative (means only) & 0.77 & 0.80 & 0.64 & 0.66 \\
			$+\ \sigma_Z$             & 0.77 & 0.78 & 0.64 & 0.66 \\
			$+\ \Delta_{ij}$          & 0.77 & 0.98 & 0.63 & 0.70 \\
			$+$ agreement             & 0.77 & 0.78 & 0.63 & 0.64 \\
			Uncertainty-aware         & 0.98 & 0.98 & 0.67 & 0.67 \\
			\bottomrule
		\end{tabular}
		\caption{Agreement of LLM mechanism labels with prescribed intervention-based reference labels across eight evidence conditions. The corresponding main-text table shows six rows. Here $n$ counts independent events, not ensemble members or LLM requests.}
		\label{tab:ladder-full-app}
	\end{table}
	
	Qualitative tags yield nearly perfect agreement in the idealized benchmark, while adding the full fingerprint is not uniformly the best condition for agreement with reference labels in the coupled ENSO experiment. In particular, GPT-4o gives 0.72 under qualitative tags and 0.67 under uncertainty-aware evidence. Because the qualitative prompt omits observation context and ensemble-derived fields supplied to quantitative prompts, such comparisons do not isolate the effect of numerical precision.
	

\subsection{One-statistic-at-a-time comparisons}
	\label{app:ablation}

The principal numerical response comparisons are accompanied by a diagnosis field distinguishing mechanism-score ambiguity from reports of inadequate observations. The fingerprint rule establishes the former diagnosis operationally. The latter was not validated by controlled changes in the observing network. Among events with no supported candidate indicator, GPT-4o reported \texttt{observations\_insufficient} for 15/15 idealized and 24/24 coupled ENSO cases; GPT-4o-mini did so for 14/15 and 22/24. Similar rates were observed under the means-only condition. These counts do not prove that the observations were objectively inadequate: a non-intervention case may have no candidate process strongly activated.

	\begin{table}[htbp]
		\centering
		\small
		\begin{tabular}{@{}llcccc@{}}
			\toprule
			& & \multicolumn{2}{c}{Reference agreement} & \multicolumn{2}{c}{Diagnosis (of 16)} \\
			\cmidrule(lr){3-4}\cmidrule(lr){5-6}
			Condition & Statistic supplied & mini & 4o & mini & 4o \\
			\midrule
			Quantitative & numeric means only & 0.77 & 0.80 & \phantom{0}1 & \phantom{0}3 \\
			$+\ \sigma_Z$ & posterior std.\ dev. & 0.77 & 0.78 & \phantom{0}1 & \phantom{0}2 \\
			$+\ \Delta_{ij}$ & pairwise separation & 0.77 & \textbf{0.98} & \phantom{0}1 & \textbf{14} \\
			$+$ agreement & ensemble agreement & 0.77 & 0.78 & \phantom{0}1 & \phantom{0}2 \\
			Uncertainty-aware & all, with verbal summary & 0.98 & 0.98 & 16 & 16 \\
			\bottomrule
		\end{tabular}
		\caption{Ablation over uncertainty components, conceptual model. Diagnosis counts ambiguous events (of $16$) identified as \texttt{mechanisms\_indistinguishable}.}
		\label{tab:ablation-app}
	\end{table}
	
	Adding score spread or ensemble agreement alone did not produce a statistically detectable improvement in either model under the supplied prompts. Supplying $\Delta_{ij}$ led to a larger change for GPT-4o than for GPT-4o-mini in the idealized cases. This does not establish that GPT-4o-mini is intrinsically unable to perform a threshold comparison; the effect of an explicit instruction to perform that arithmetic was not tested.

\subsection{Exploratory response characteristics}
	\label{app:null}
	
	The exploratory response analysis considered whether the model named the competing mechanisms, referred to observation quality, quoted a numerical statistic, matched expressed confidence to ensemble agreement, or suggested additional observations. The reported differences were not stable across conditions. Suggestions for additional measurements especially depend on which state variables the prompt explicitly identifies as observed; they do not independently validate an observing-system recommendation. We do not use these exploratory categories as formal performance metrics.

\section{Differences between the controlled benchmarks}
	\label{app:tier1}

    On the conceptual model, gated mechanisms exceed the threshold by a median of $3.4$ (Appendix~\ref{app:c1_specific}) posterior standard deviations (range $2.0$ to $9.0$), so support is comparatively clear in this benchmark; qualitative tags are often sufficient for the tested attribution task, although they do not carry the full joint-uncertainty information. On the ICM the margin is one to three standard deviations. A test that thinned the conceptual model's observations (interval $20\to80$) to reproduce the ICM regime did not work: gated $\bar Z$ stayed at $6$--$11$ in this sensitivity experiment. The observation operators and mechanisms differ between the two benchmarks, so changing one observation interval need not make their score distributions comparable. Reducing the imposed amplitude would be another possible sensitivity experiment, but it was not tested. The difference in LLM performance cannot be assigned uniquely to dynamical complexity because event composition, perturbations, observations, and thresholds also differ.
	

\section{Illustrated coupled ENSO events and their fingerprints}
\label{app:uqevents}

This section describes how the ICM events in Figure~\ref{fig:saicm-diagnostics-app} were selected and examines the corresponding posterior-uncertainty figures (Figures~\ref{fig:saicm-diagnostics-app}b and~\ref{fig:uq_blockc}). It also details how the fingerprint rule can give a different label from a point estimate.

\subsection{How the events were selected}
\label{app:uq_selection}

\paragraph{Event library.} The ICM library contains 90 independent twin simulations covering five event types: thermocline-driven, wind-burst-driven, advection-driven, mixed thermocline–advection, and null events, with various realizations of each type. Each realization begins with a background period without imposed process amplification, followed by an attribution window in which selected processes are amplified. The imposed intervention defines the reference label; the SST response illustrates its effects but is not used as a threshold for assigning that label. For each realization, the assimilation system stores the full posterior ensemble ($M=150$ members) for SST, the Kelvin and Rossby amplitudes, and the wind-burst amplitude, together with the truth run. The fingerprint ($\bar Z_k$, $\sigma_{Z_k}$, $\Delta_{ij}$) and its label are computed from the posterior ensemble. The reference simulations and paired counterfactuals are used to examine responses and to display physical fields; intervention labels are set by the activated gates.

\paragraph{Panel (a) of Figure~\ref{fig:saicm-diagnostics-app}: a representative correctly labeled event per process.} The panel illustrates selected cases in which the fingerprint agrees with the intervention-based reference label; it is not an unbiased sample of all realizations. For each prescribed mechanism category we therefore took the realizations for which the fingerprint label equals the prescribed intervention label and which have a stored posterior, and chose the one whose margin is closest to the median margin of those realizations. The margin is
\begin{equation}
 m \;=\;
 \begin{cases}
 \bar Z_{\rm true}-\max_{k\neq {\rm true}}\bar Z_k & \text{single-mechanism events},\\[2pt]
 -\,|\bar Z_{\rm thermocline}-\bar Z_{\rm advective}| & \text{mixed events (a smaller gap is a cleaner mixed event).}
 \end{cases}
 \label{eq:margin}
\end{equation}

\paragraph{Figure~\ref{fig:uq_blocka}: a mixed event and a single-mechanism contrast.}  LLMs struggle to attribute the two events without uncertainty (thermocline-driven, realization~4; mixed thermocline--advective, realization~5). The mixed event is one in which the larger mean score (the point-estimate answer) names only one of the two mechanisms, so that the separation $\Delta$ is what makes the label mixed. The thermocline-driven event is the contrast: a single mechanism is supported, so $\Delta$ plays no role in its label. 

\subsection{The uncertainty's role}
\label{app:uq_reading}

Figures~\ref{fig:saicm-diagnostics-app}b and~\ref{fig:uq_blockc} show the posterior standard deviation of each plotted field. The standard deviation is taken across the $M=150$ members in physical units. The right-hand column gives the fingerprint of the event. Bars are $\bar Z_k$ with $\pm\sigma_{Z_k}$ error bars, the dashed line is the support threshold $z_{\rm thr}=1$, and the bracket gives $|\Delta_{ij}|$ between the two mechanisms with the highest $\bar Z$. Only a pair in which \emph{both} mechanisms exceed the threshold enters the label (Section~\ref{sec:framework}); when the runner-up is below the threshold, the bracket is shown for completeness and plays no role in the decision. 

Three features are common to all rows. First, the standard deviation (SD) of $H$ and $U$ is concentrated in the west, roughly west of $200^\circ$E, with the largest values near the western boundary and in the warm pool. This pattern is consistent with weaker SST–thermocline coupling in the west than in the east. The spread maps alone do not determine whether this pattern arises from spatial differences in observability, forcing uncertainty, or the coupled dynamics. Second, the SD maps exhibit horizontal banding with a period of about one month, consistent with the monthly analysis cycle. The spread also tends to be larger where the signal is stronger. Accordingly, the events with the strongest thermocline or wind signals (rows 1 and 2 of Figure \ref{fig:saicm-diagnostics-app}) have the largest $H$ and wind SDs. The advective event (row 3), whose $H$ signal is weak, has the smallest spread.

\begin{figure}[p]
\centering
\includegraphics[width=\textwidth]{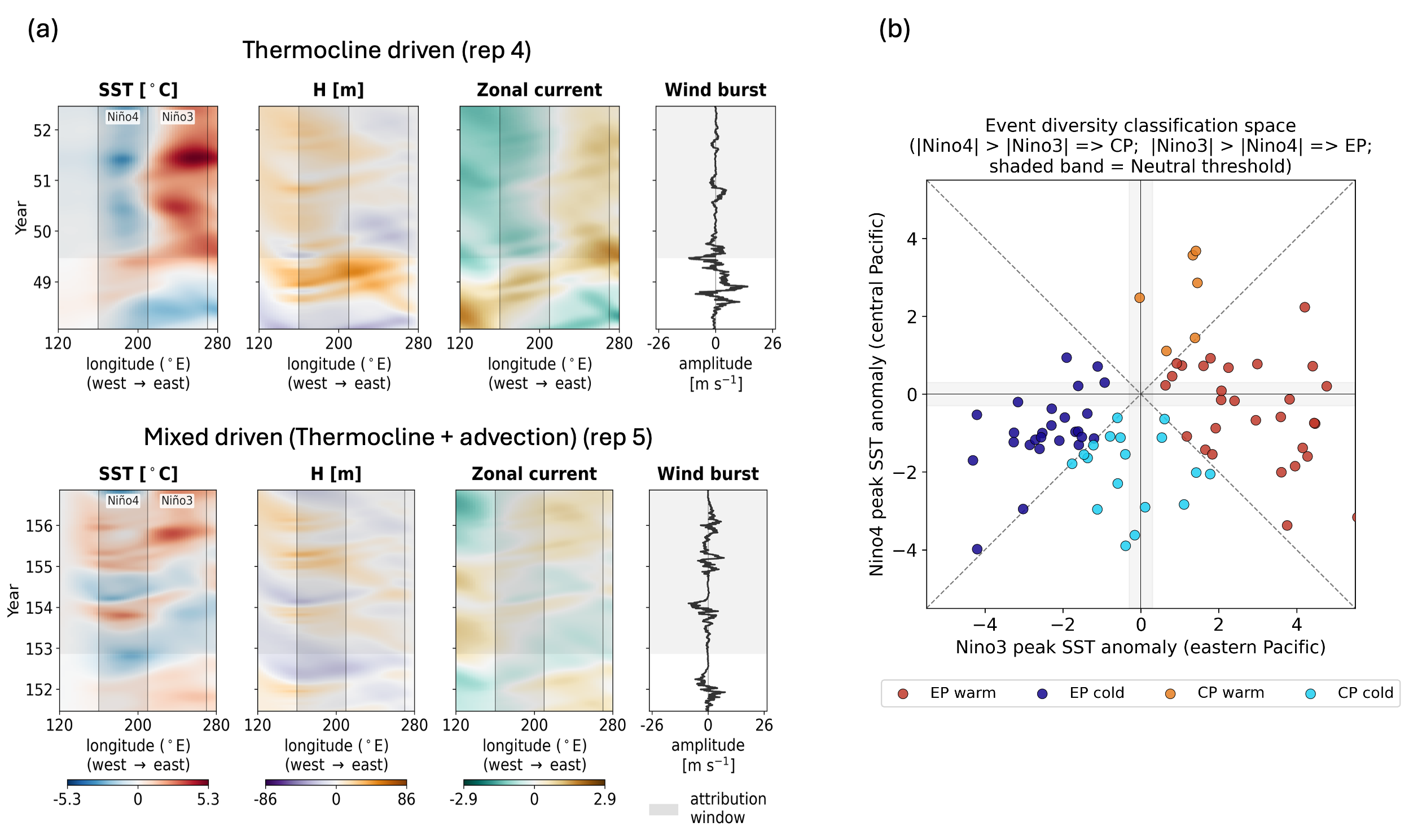}
\caption{(a) A thermocline-reference event (top) and a mixed thermocline--advection reference event (bottom), with different fingerprint labels. The shaded longitude band denotes Ni\~{n}o-4, and vertical lines bound Ni\~{n}o-3. (b) Peak Ni\~{n}o-3 and Ni\~{n}o-4 SST anomalies show the diversity of simulated SST patterns; diagonal boundaries distinguish eastern- and central-Pacific classifications. }
\label{fig:uq_blocka}
\end{figure}

\begin{figure}[p]
\centering
\begin{subfigure}{0.93\textwidth}
\includegraphics[width=\textwidth]{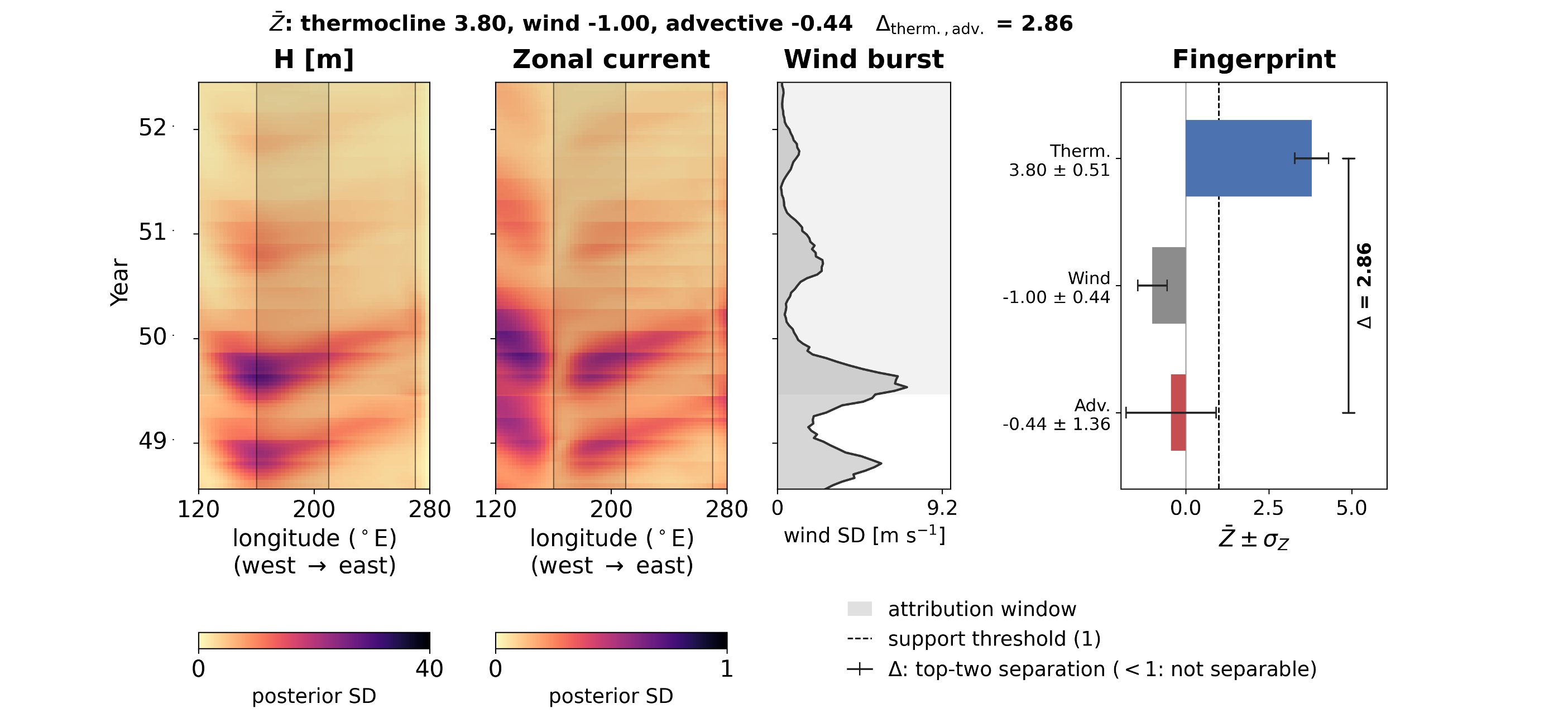}
\caption{Thermocline-driven event of Figure~\ref{fig:uq_blocka}a (realization 4).}
\end{subfigure}\\[6pt]
\begin{subfigure}{0.93\textwidth}
\includegraphics[width=\textwidth]{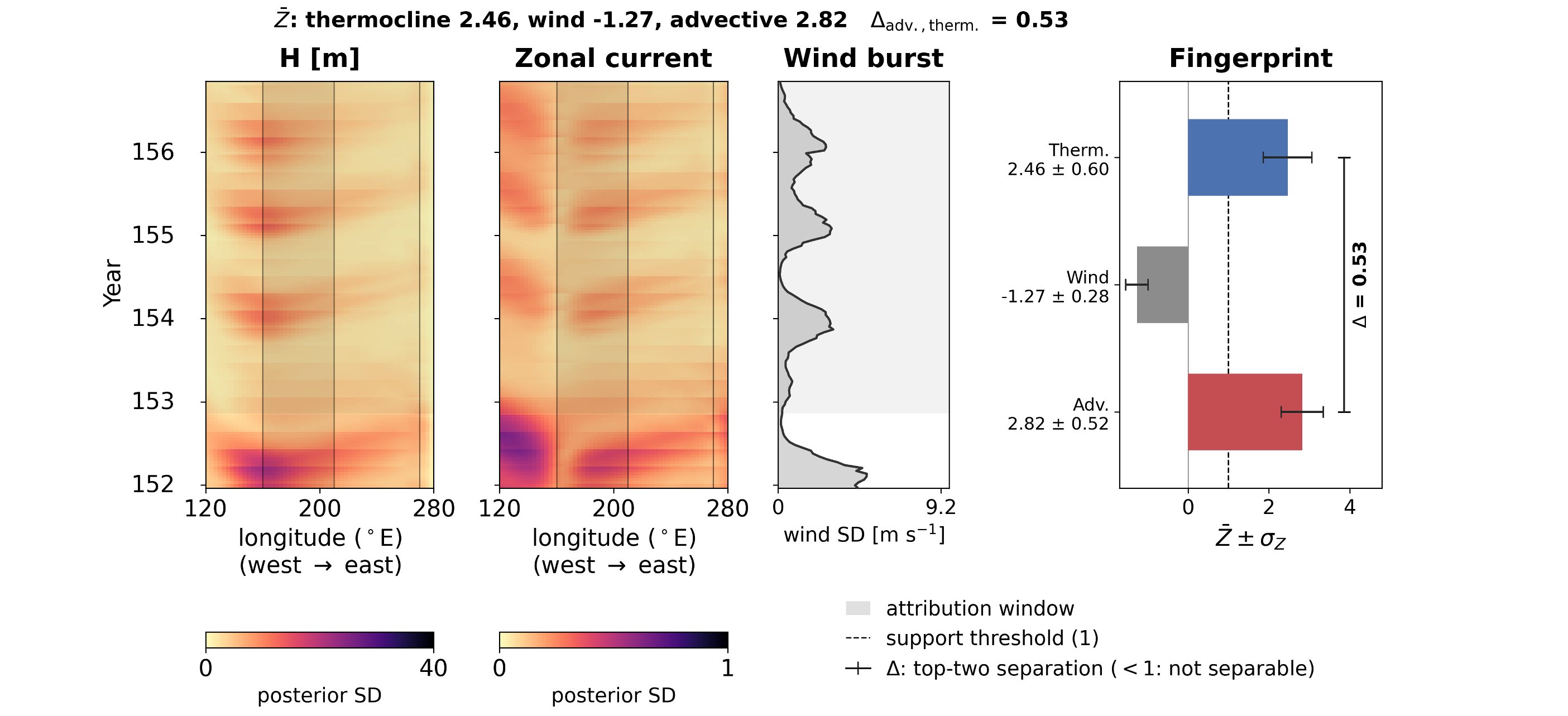}
\caption{Mixed-driven event (thermocline--advection) of Figure~\ref{fig:uq_blocka}a (realization 5).}
\end{subfigure}
\caption{Posterior standard deviation across the 150 ensemble members of the two events of Figure~\ref{fig:uq_blocka}a. The uncertainty is in physical units and on separate color scales for each variable. Indicator scores are normalized and dimensionless, making them comparable across mechanisms.}
\label{fig:uq_blockc}
\end{figure}

Additional examples that could not be attributed without providing uncertainty are shown in Figure~\ref{fig:uq_blocka}a, together with their standard deviations and indicator scores in Figure~\ref{fig:uq_blockc}. Figure~\ref{fig:uq_blocka}b illustrates the diversity of SST patterns in the Ni\~{n}o-3--Ni\~{n}o-4 classification space used in the tests.

\subsection{Why each event receives its label}
\label{app:uq_decisions}

Table~\ref{tab:uq_decision} lists the indicators for events in Figures~\ref{fig:saicm-diagnostics-app} and ~\ref{fig:uq_blocka}. The label follows the rule of Section~\ref{sec:framework}: a mechanism is supported if $\bar Z_k>1$; if two supported mechanisms have $|\Delta_{ij}|<1$ the label is mixed; otherwise, it is the highest supported mechanism. The column $P(Z_k>1)$ is the probability of exceeding the threshold under a Gaussian approximation, $\Phi\big((\bar Z_k-1)/\sigma_{Z_k}\big)$; it is a summary of the bars, not a separate statistic of the ensemble.

\begin{table}[htbp]
\centering
\small
\setlength{\tabcolsep}{4pt}
\caption{Fingerprints and labels of the events driven by different mechanisms in figures~\ref{fig:saicm-diagnostics-app} and~\ref{fig:uq_blocka}. Th, W, Ad: thermocline, wind burst, advective. $\bar Z_k\pm\sigma_{Z_k}$ and $P(Z_k>1)$ (Gaussian approximation) for each mechanism. $\Delta$: separation of the two highest $\bar Z$ (it enters the label only if both are supported). Fingerprint: label from the fingerprint rule (Section~\ref{sec:framework}); point est.: the deterministic classifier, i.e.\ the highest $\bar Z$. $\checkmark$/$\times$: agreement with the prescribed intervention label.}
\label{tab:uq_decision}
\resizebox{\textwidth}{!}{%
\begin{tabular}{@{}lcccccccl@{}}
\toprule
Driven by  & $\bar Z_{\rm Th}$ & $\bar Z_{\rm W}$ & $\bar Z_{\rm Ad}$ & $P(Z_k>1)$ Th/W/Ad & $|\Delta|$ & Supported & Fingerprint & Point est. \\
\midrule
\multicolumn{9}{@{}l}{\emph{Panel a of Figure~\ref{fig:saicm-diagnostics-app}}}\\
Thermocline  & $5.00\pm1.01$ & $2.28\pm0.87$ & $-2.31\pm1.29$ & 1.00/0.93/0.01 & 2.79 & Th, W & Th $\checkmark$ & Th $\checkmark$ \\
Wind burst  & $-0.98\pm0.43$ & $2.19\pm0.98$ & $-1.06\pm0.58$ & 0.00/0.89/0.00 & 3.64 & W & W $\checkmark$ & W $\checkmark$ \\
Advective  & $-4.92\pm1.18$ & $-0.89\pm0.36$ & $4.33\pm1.30$ & 0.00/0.00/0.99 & 3.93 & Ad & Ad $\checkmark$ & Ad $\checkmark$ \\
Mixed  & $2.02\pm0.51$ & $-0.72\pm0.42$ & $1.64\pm2.39$ & 0.98/0.00/0.61 & 0.16 & Th, Ad & mixed $\checkmark$ & Th $\times$ \\
\addlinespace
\multicolumn{9}{@{}l}{\emph{Figure~\ref{fig:uq_blocka} }}\\
Thermocline  & $3.80\pm0.51$ & $-1.00\pm0.44$ & $-0.44\pm1.36$ & 1.00/0.00/0.14 & 2.86 & Th & Th $\checkmark$ & Th $\checkmark$ \\
Mixed  & $2.46\pm0.60$ & $-1.27\pm0.28$ & $2.82\pm0.52$ & 0.99/0.00/1.00 & 0.53 & Th, Ad & mixed $\checkmark$ & Ad $\times$ \\
\bottomrule
\end{tabular}}
\end{table}

\paragraph{Thermocline-driven event (Figure~\ref{fig:saicm-diagnostics-app}a, row 1).} The thermocline score is well above the threshold ($5.00\pm1.01$, $P\approx1$). The wind score is also supported ($2.28\pm0.87$, $P\approx0.93$), so both mechanisms initially appear plausible. However, the separation is $\Delta=2.79$, well above the threshold of $1$, and the event is therefore labelled thermocline-driven. Although the two elevated means could suggest a mixed event, the separation indicates that the thermocline indicator ranks well above the wind indicator under the adopted rule. The advective score is $2.31$ background standard deviations \emph{below} the background mean (Section~\ref{app:uq_roles}), providing evidence against advection.

\paragraph{Wind-burst-driven event (row 2).} Only the wind score is supported ($2.19\pm0.98$); the thermocline and advective scores are negative. The event is therefore labelled wind-burst-driven, without using $\Delta$; the value $3.64$ in brackets gives the gap to the runner-up. This is the least secure single-mechanism attribution among the four events. The relative spread, $\sigma_Z/\bar Z=0.45$, is larger than the corresponding values for the dominant mechanisms in the thermocline and advective events ($0.2$--$0.3$). Moreover, $P(Z_{\rm wind}>1)\approx0.89$, corresponding to an approximate probability of 0.11 that the wind score falls below the threshold under the assumed Gaussian distribution. The large wind SD in this row, which approaches the axis maximum, is consistent with this uncertainty.

\paragraph{Advective event (row 3).} The advective score is the only supported score ($4.33\pm1.30$, $P\approx0.99$). The thermocline score is $4.92$ background standard deviations below the background mean. The event is therefore labelled advective.

\paragraph{Mixed event (row 4).} Both the thermocline and advective scores are supported: $2.02\pm0.51$ and $1.64\pm2.39$, respectively. Because their separation is $\Delta=0.16<1$, the event is labelled mixed. A point estimate would instead classify it as thermocline-driven, since $2.02>1.64$.

Two details are important. First, the advective score is weakly constrained: its error bar spans $[-0.75,4.03]$, and $P(Z_{\rm Ad}>1)\approx0.61$. Second, the large uncertainty is directly relevant to the standardized separation. For any correlation, $\mathrm{Var}(Z_i-Z_j)\ge(\sigma_{Z_i}-\sigma_{Z_j})^2$, $|\Delta_{ij}|\le|\bar Z_i-\bar Z_j|/|\sigma_{Z_i}-\sigma_{Z_j}|=0.38/1.88=0.20$. The large advective spread itself contributes substantially to this small separation, and no admissible correlation could increase $|\Delta|$ above $0.20$. The resulting mixed fingerprint label reflects limited separation under the prescribed rule; it does not establish that both physical processes contribute equally to the event.

\paragraph{Events of figure~\ref{fig:uq_blockc}.} For the thermocline-driven event, only the thermocline score is supported ($3.80\pm0.51$), while the advective score is negative ($-0.44\pm1.36$). The event is therefore labelled thermocline-driven. The separation is not used; the value $2.86$ in brackets is the gap to the runner-up.

The mixed event provides the clearer of the two mixed examples. Both scores are well supported ($2.46\pm0.60$ and $2.82\pm0.52$, with $P\ge0.99$), and their difference, $0.36$, is smaller than either spread. Consequently, $\Delta=0.53<1$, and the event is labelled mixed.
\subsection{What the uncertainty changes}
\label{app:uq_roles}

\paragraph{From rank to separation.} The two mixed events (Figure~\ref{fig:saicm-diagnostics-app}, row 4, and Figure~\ref{fig:uq_blockc}b) have the same prescribed mechanisms. However, the thermocline score is larger in the first event, whereas the advective score is larger in the second. In both cases, the difference between the scores ($0.38$ and $0.36$) is smaller than their respective spreads ($0.51$ and $2.39$; $0.60$ and $0.52$). The leading mean score can therefore differ across stochastic event realizations with the same prescribed intervention class. 

A point-estimate attribution therefore gives different labels for two events with the same prescribed intervention class, while the fingerprint-based method labels both as mixed. For comparison, the thermocline event in row 1 also has two supported scores. Their separation is sufficiently large to give a single-process fingerprint label, although both processes can still influence the physical evolution.

\paragraph{Support is probabilistic.} A score mean above $1$ meets the formal support criterion, but the corresponding evidence can be weak if its posterior spread is large. For the scores used in Table~\ref{tab:uq_decision}, $P(Z_k>1)$ ranges from approximately $0.61$ for the supported mechanisms of the mixed event in row 4 to nearly $1.00$. Because the label is based on the mean scores, it is unchanged when a small number of ensemble members fall below the threshold. However, the strength of the evidence varies across labels. A single-mechanism label depends on one score, with $P\ge0.89$ in these examples, whereas a mixed label requires both scores to exceed the threshold and is therefore limited by the weaker score, which is $0.61$ in row 4.

\paragraph{Negative scores.} A negative score, $\bar Z_k<0$, means that an indicator is below its background reference level in the reconstructed ensemble. This provides evidence against unusually elevated activity under that indicator, not proof that the physical process is absent. The criterion is one-sided: a mechanism that is weaker than the background is not considered a candidate. For this reason, the decision uses $\bar Z_k$ rather than $|\bar Z_k|$. Taking the absolute value would incorrectly treat the advective score of the thermocline event ($-2.31$) as evidence for advection. In the advective event, it would also make the thermocline score ($|-4.92|$) larger than the advective score ($4.33$), leading to an incorrect thermocline attribution.

\paragraph{Spread scales with the score.} For the dominant mechanism in each event, the ratio $\sigma_Z/\bar Z$ ranges from $0.13$ to $0.30$, except for the wind-burst event, where it is $0.45$. The spread of each score increases approximately with its magnitude, as expected for an indicator based on the integral of the absolute value of the reconstructed field. The support rule in this paper uses the mean score and therefore does not itself adjust the support decision for the spread. A larger score is therefore not necessarily more certain; the relevant quantity is the standardized distance from the threshold $(\bar Z_k-z_{\rm thr})/\sigma_{Z_k}$.

\subsection{An example of a false ambiguity flag}
	\label{app:sacrifice}
	
On the thermocline-driven event of Figure~\ref{fig:saicm-diagnostics-app}a (row~1; ICM, GPT-4o) the prescribed reference label is thermocline and the event is unambiguous by the pipeline's criterion: the thermocline and wind-burst scores are both supported ($\bar Z_{\rm thermocline}=5.00\pm1.01$, $\bar Z_{\rm wind\_burst}=2.28\pm0.87$), but their separation, $\Delta=2.79$, exceeds $\delta_{\rm thr}=1.0$ (Table~\ref{tab:uq_decision}). Given qualitative tags, the model returns \texttt{mixed} (``Both thermocline and wind burst contributions are high $\ldots$ the evidence does not clearly distinguish''). Given the uncertainty-aware fingerprint, the same model returns \texttt{thermocline} (``clearly separated from the wind\_burst mechanism''). Both prompts describe the same simulated event, but they differ in more than the separation statistic: the uncertainty-aware prompt also contains additional ensemble information and a verbal conclusion. This example illustrates different outputs, but does not isolate the effect of any single prompt field.

\section{Additional controlled ENSO events}
\label{app:additional-enso}

\begin{figure}[p]
\centering
\includegraphics[width=\textwidth]{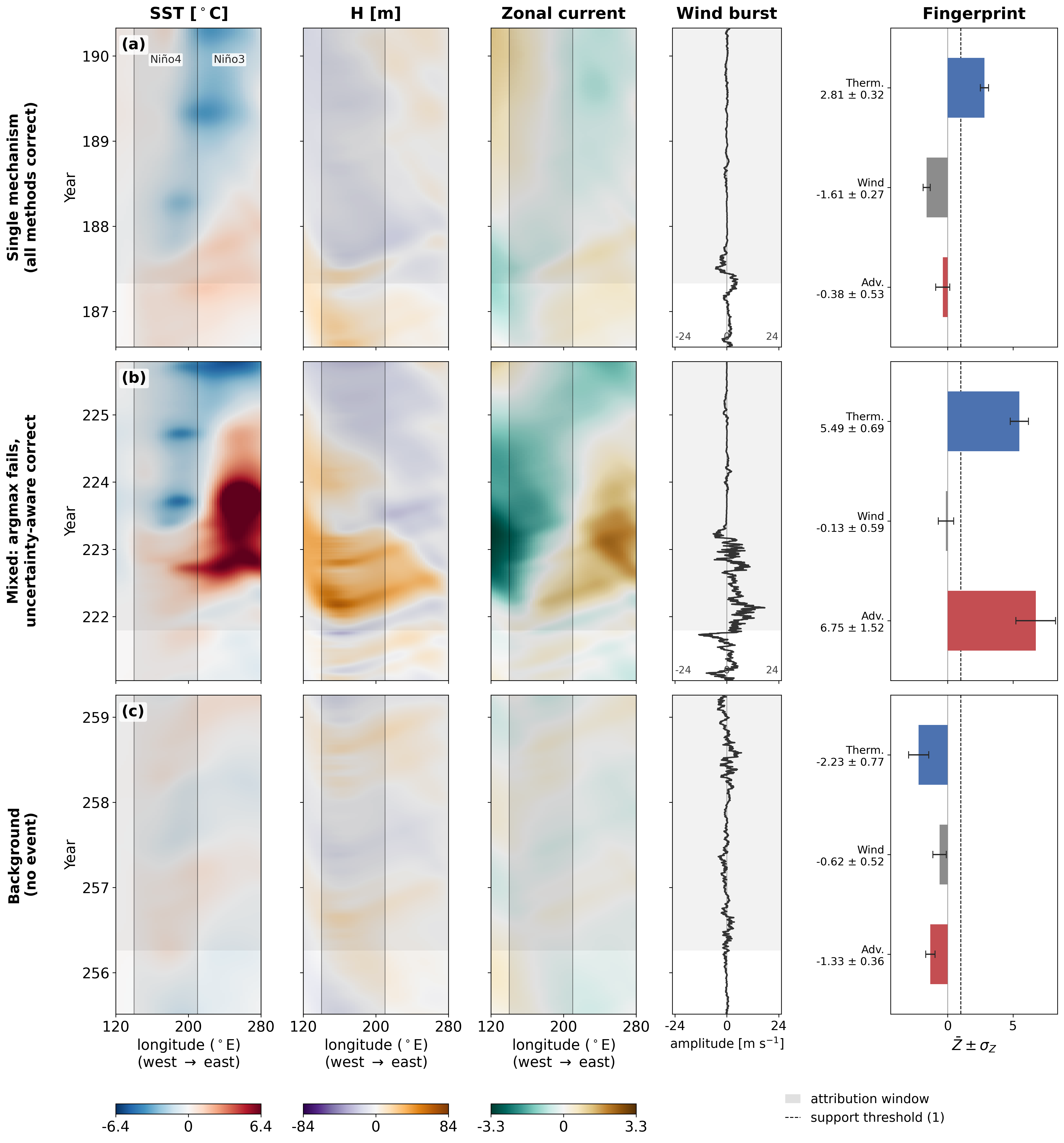}
\caption{Additional ICM events, rows (a)--(c) of Table~\ref{tab:icm_extra}. Columns: SST, thermocline depth $H$, zonal current $U$ (Hovm\"oller diagrams against longitude) and the wind-burst amplitude (time series). The vertical axis is years from event onset. The shaded band is the attribution window. The shaded longitude band is the Ni\~no-4 region, and the vertical lines bound the Ni\~no-3 region. Colour scales and the wind-amplitude axis are shared across all rows and both pages. Fields are the simulated (truth) fields of the library run.}
\label{fig:icm_hov_1}
\end{figure}

\begin{figure}[p]
\centering
\includegraphics[width=\textwidth]{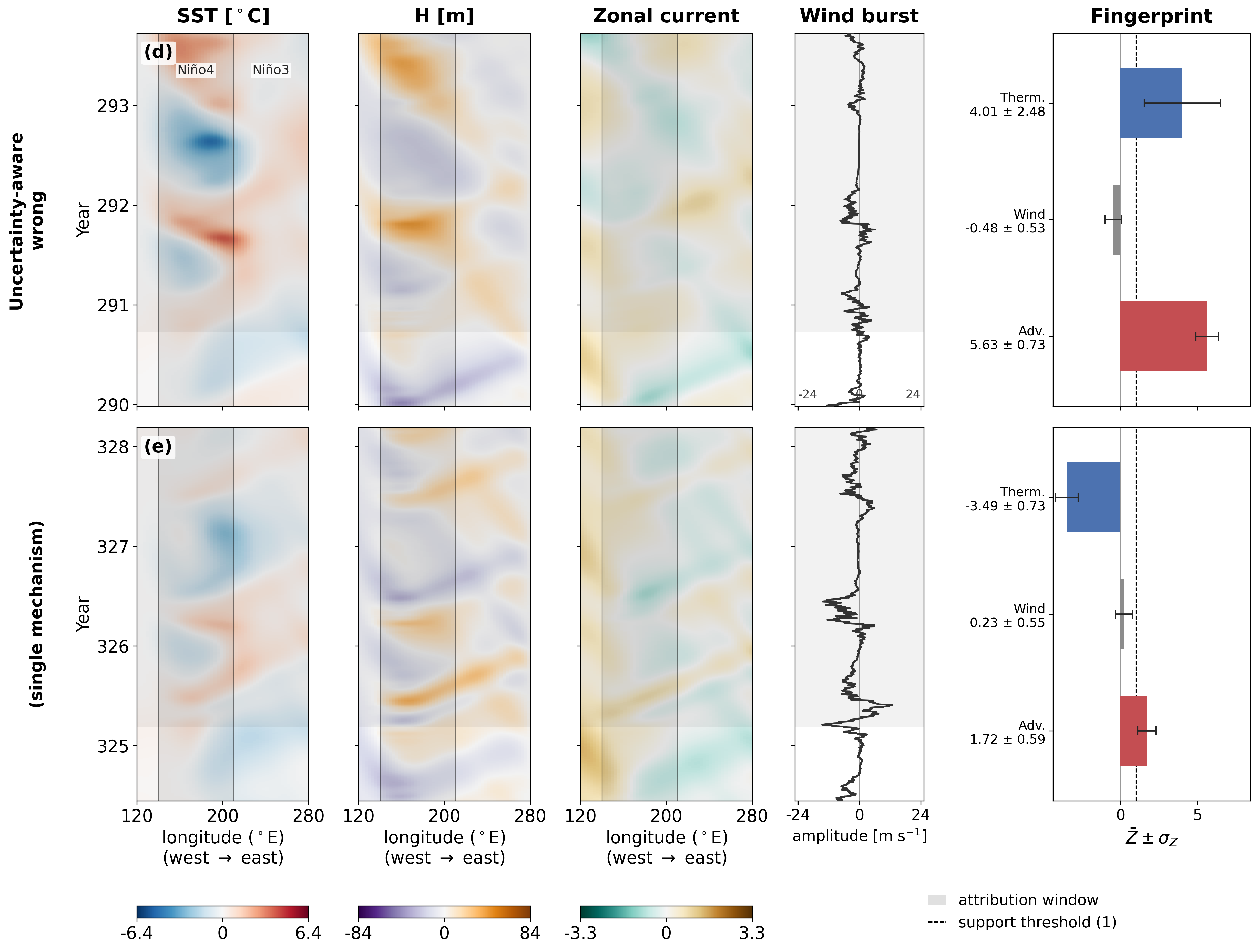}
\caption{As Figure~\ref{fig:icm_hov_1}, for rows (d) and (e) of
Table~\ref{tab:icm_extra}.}
\label{fig:icm_hov_2}
\end{figure}

\begin{figure}[t]
\centering
\includegraphics[width=0.7\textwidth]{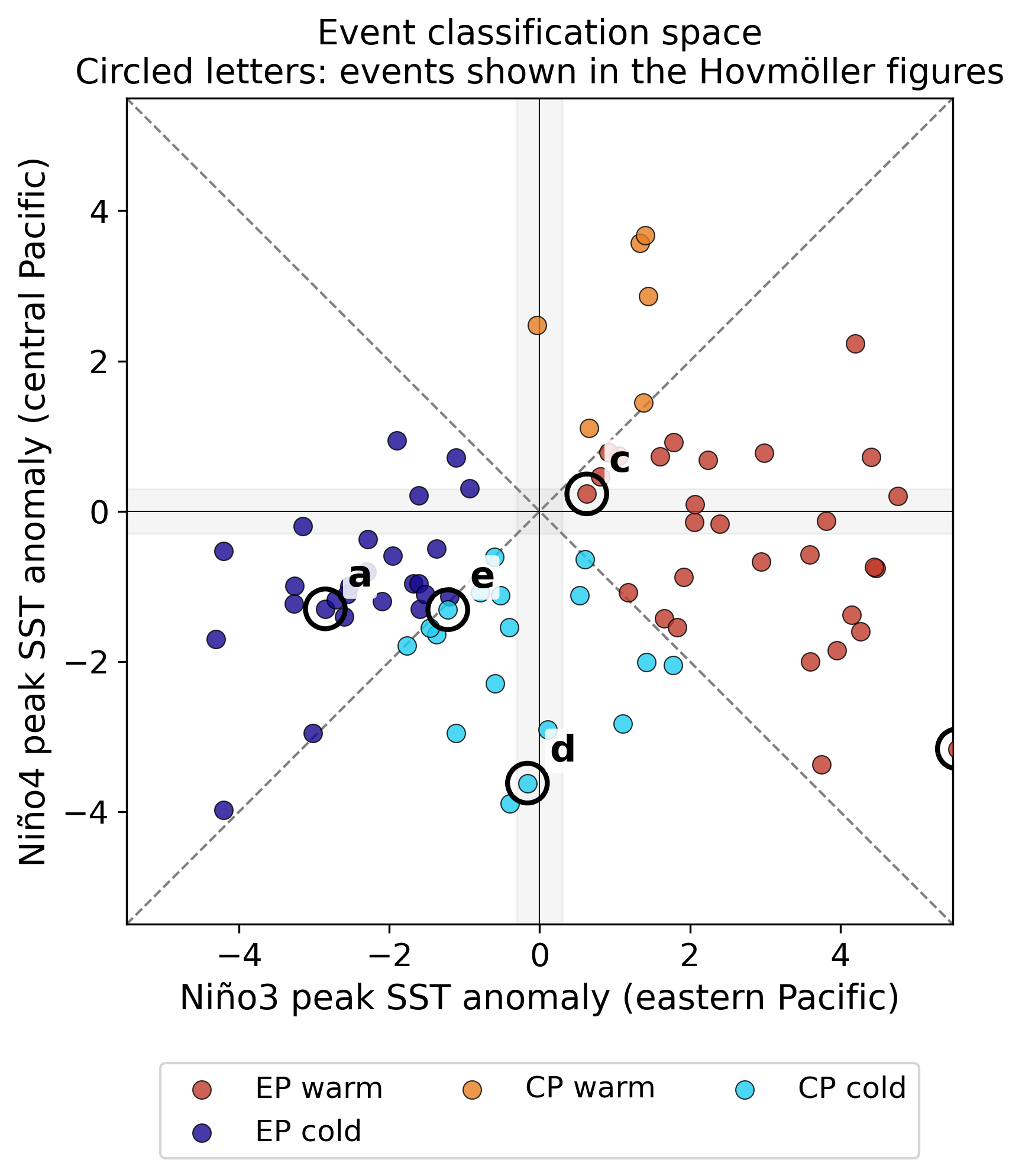}
\caption{Ni\~no-3 and Ni\~no-4 peak SST anomalies for all events in the ICM library, coloured by SST class, with the events of Table~\ref{tab:icm_extra} circled and lettered as in Figures~\ref{fig:icm_hov_1} and~\ref{fig:icm_hov_2}. The dashed lines are $|\text{Ni\~no 3}|=|\text{Ni\~no-4}|$ and the shaded bands mark the neutral threshold.}
\label{fig:icm_scatter}
\end{figure}

Figures~\ref{fig:icm_hov_1} and~\ref{fig:icm_hov_2} show five further events from the ICM twin-experiment library, and Figure~\ref{fig:icm_scatter} shows where they sit in the library's Ni\~no-3/Ni\~no-4 classification space. The events were chosen by a deterministic selection rule applied to the stored fingerprint table. In each class of Table~\ref{tab:icm_extra} we took the first event, in ascending realization index, that met the criterion, excluding the events shown in the main text. Within each event, the model's own state is plotted as stored: SST, thermocline depth $H=K-R$, zonal current $U=K+R$, and the wind-burst amplitude as a time series.

The selection includes cases in which the fingerprint agrees with intervention labels, as well as cases illustrating a mismatch or unresolved attribution. Event (a) is a single-process intervention for which the reported fingerprint agrees with the reference label. Event (b) is a mixed event in which the argmax attribution is incorrect, whereas the fingerprint-based attribution agrees with the intervention label. Event (c) is a background control with no deliberately activated event. Event (d) is a case in which the uncertainty-aware label disagrees with the prescribed perturbation-based reference label. Event (e) is another single-process intervention with a cold SST peak.

Four of these selected realizations have negative SST anomalies. The sign depends on the particular stochastic realization and is not fixed by the prescribed process alone (Table~\ref{tab:icm_extra}).

Two features of the table are important for interpreting the figures. First, the SST class and the prescribed intervention label are separate labels. For example, the thermocline-only event (a) exceeds the neutral threshold in both the Niño-3 and Niño-4 regions and is therefore classified as mixed CP--EP. Thus, a single perturbed process can generate large SST anomalies in both regions; this SST pattern label is distinct from the mixed mechanism-fingerprint label.

Second, event (b) has Niño-3 and Niño-4 anomalies with opposite signs ($+5.56$ and $-3.16$\,K). It is therefore characterized by eastern warming and central cooling, rather than by an El Ni\~{n}o pattern with warming in both regions. The ``mixed CP--EP'' label indicates only that both regions exceed the threshold; it does not specify the sign of the anomalies. Event (b) also has the largest anomalies of the five events. Because the colour scales are shared across all rows, the weaker events appear faint by comparison. Event (c) has the smallest anomalies ($-0.54$ and $+0.04$\,K), so its eastern classification is based on a relatively weak signal.

\begin{table}[t]
\caption{Additional ICM events shown in Figures~\ref{fig:icm_hov_1} and~\ref{fig:icm_hov_2}. Sign and SST class follow the convention used for the library (peak of the basin-mean absolute SST within the event window; class from the Ni\~no-3 and Ni\~no-4 anomalies relative to the neutral threshold). Anomalies are in $^\circ$C.}
\label{tab:icm_extra}
\centering
\begin{tabular}{lllccrr}
\hline
Row & Prescribed intervention & Realization & Sign & SST class & Ni\~no-3 & Ni\~no-4 \\
\hline
(a) & thermocline only        & 0 & cold & mixed CP--EP & $-2.85$ & $-1.30$ \\
(b) & thermocline + advective & 7 & warm & mixed CP--EP & $+5.56$ & $-3.16$ \\
(c) & background (no event)   & 0 & cold & EP           & $-0.54$ & $+0.04$ \\
(d) & advective only          & 0 & cold & CP           & $-0.16$ & $-3.62$ \\
(e) & advective only          & 1 & cold & mixed CP--EP & $-1.23$ & $-1.31$ \\
\hline
\end{tabular}
\end{table}

\section{An illustrated LLM evidence-ladder example}
	\label{app:trace}

\texttt{mixed\_thermo\_wind\_0} is an archived identifier for an idealized event in which $Y_1$ and $Y_2$ are perturbed together. The mnemonic names in the identifier and original prompts do not imply that these drivers are physical thermocline or wind processes. With the observation-only input, neither model receives an explicit mechanism fingerprint. They also defer to the point estimate and differ only in their confidence: 4o reports low confidence and \texttt{EVIDENCE\_SUFFICIENT: no}, whereas mini reports high confidence. Both interpret the qualitative tags as \emph{mixed}.

Without the verbal summary, mini describes the separation as ``adequate'' and identifies a single mechanism with confidence, despite $\Delta_{\rm thermocline,wind}=1.04$ being below the stated threshold $\delta_{\rm thr}=2.5$. In contrast, 4o compares $1.04$ with $2.5$ and classifies the result as \emph{mixed}. When the verbal summary is provided, both models converge on \emph{mixed} and cite $\Delta=1.04$ in their explanations. The archive contains the full prompts and responses for all eight conditions and both models.

\section{Historical ENSO application using GODAS}
\label{app:godas_events}
\subsection{GODAS assimilation and score normalization}
In the historical application, monthly GODAS SST from 1981--2020 is assimilated into the intermediate coupled ENSO model. The reported configuration uses 150 members, RTPS relaxation, and SST-observation error standard deviation $0.3$ K. Fingerprints are computed from the assimilated ensemble over a fixed 12-month window (March of the reported peak year through the following February). Such temporal aggregation represents processes evolving around an SST event; it does not by itself identify an optimal physical lead--lag relation.

The observation operator maps model SST to anomalies relative to the model's own monthly climatology, computed from a 60-year free run, and these are compared with GODAS monthly SST anomalies averaged over $5^\circ$S--$5^\circ$N and interpolated linearly in longitude onto the model grid. Only grid points within the longitudinal coverage of the GODAS record ($140^\circ$E--$280^\circ$E) are observed. Observations are assimilated monthly. No mechanism gates are prescribed in the historical run. The regression-based thermocline and advective indicators therefore measure effective departures of the reconstructed SST tendency from baseline model feedbacks, rather than known changes in physical feedback coefficients. These residuals may also reflect model mismatch or unresolved forcing. Their scores are thus not directly comparable in magnitude with the controlled gate experiments.

For this historical application, the reference means and standard deviations used to standardize mechanism scores are computed from all 12-month windows of the 40-year assimilated record. The thresholds from the controlled ENSO benchmark are applied without retuning. Unlike the controlled benchmark, this reference distribution includes periods with observed ENSO variability, so the standardized scores and the transferred thresholds need not have the same interpretation as in the controlled experiments.

\subsection{Historical event selection}

Table~\ref{tab:godas_events} lists the eight events analyzed in the observed-SST application. They were selected from GODAS SST anomalies for 1981--2020 by a fixed rule: an event is a period with a three-month mean Ni\~no-3.4 anomaly beyond $\pm0.5$\,K for at least five consecutive months, and its flavor (eastern or central) follows the Ni\~no-3 and Ni\~no-4 combination of \citet{renjin2011} at the peak. The rule gives 20 events; two events were selected in each of four classes according to the SST flavor-margin criterion, namely eastern and central El~Ni\~no and eastern and central La~Ni\~na. The reference mechanism for this LLM-fidelity comparison is the fingerprint label (column Fingerprint), computed from the assimilated ensemble with the rule of Section~\ref{sec:framework}; it is not inferred directly from the SST flavor of the event and is not an independently verified causal attribution. The table also gives whether the uncertainty-aware condition and the other evidence conditions reproduce it. The classification uses the \citet{renjin2011} indices $N_{CT}=N_3-\alpha N_4$ and $N_{WP}=N_4-\alpha N_3$, where $\alpha=0.4$ if $N_3N_4>0$ and $\alpha=0$ otherwise. The reported signed flavor margin is $F=|N_{CT}|-|N_{WP}|$. For central-Pacific events, this signed margin is commonly negative; the ranking rule therefore needs to account for its sign. The 1994 event has a small absolute SST-flavor margin, so it lies relatively close to the CP/EP boundary; this is distinct from ambiguity in its physical-mechanism attribution.

\begin{table}[t]
\caption{Events used in the observed-SST application, selected from GODAS 1981--2020 as two events per class under the stated flavor-margin ranking (see selection details above). EP and CP denote eastern- and central-Pacific events. Ni\~no-3 and Ni\~no-4 are SST anomalies (K) averaged over the three months centred on the peak of the Ni\~no-3.4 index; for these events the class coincides with the larger of $|$Ni\~no-3$|$ and $|$Ni\~no-4$|$. \emph{Fingerprint}: the label of the pipeline (Section~\ref{sec:framework}), computed from the assimilated ensemble, which serves as the reference. \emph{UA}: uncertainty-aware condition. The \emph{Other} column is an aggregate from the original experiment record that requires a precise model-by-condition definition. A yes value signifies agreement with the model-based fingerprint reference, not independent verification of the mechanism. 
\emph{UA}: uncertainty-aware condition (full fingerprint with verbal summary). \emph{Other}: Point-wise or the $+\Delta_{ij}$ condition (numerical fingerprint with the pairwise separation but without the verbal summary). }
\label{tab:godas_events}
\centering
	\begin{tabular}{llrrlcc}
			\hline
			Event extent & Class & Ni\~no-3 & Ni\~no-4 & Fingerprint  & Other & UA \\
			\hline
			1994-09 to 1995-02 & CP El~Ni\~no & $+0.86$ & $+0.96$ & advective    & no& no  \\
			1995-08 to 1996-04 & EP La~Ni\~na & $-1.11$ & $-0.67$ & thermocline  & no  & yes \\
			1997-05 to 1998-04 & EP El~Ni\~no & $+3.62$ & $+0.74$ &  thermocline + wind & no  & yes \\
			1998-06 to 2000-06 & CP La~Ni\~na & $-1.31$ & $-1.61$ &  thermocline & yes & yes \\
			2004-08 to 2005-01 & CP El~Ni\~no & $+0.44$ & $+0.85$ & advective   & no  & yes \\
			2011-08 to 2012-02 & CP La~Ni\~na & $-0.89$ & $-1.13$ & advective    & no  & yes \\
			2014-10 to 2016-04 & EP El~Ni\~no & $+2.88$ & $+1.45$ & wind & no  & yes \\
			2017-10 to 2018-03 & EP La~Ni\~na & $-1.05$ & $-0.50$ & thermocline  & no  & yes \\
			\hline
	\end{tabular}
\end{table}

\begin{figure}[t]
\centering
\includegraphics[width=\textwidth]{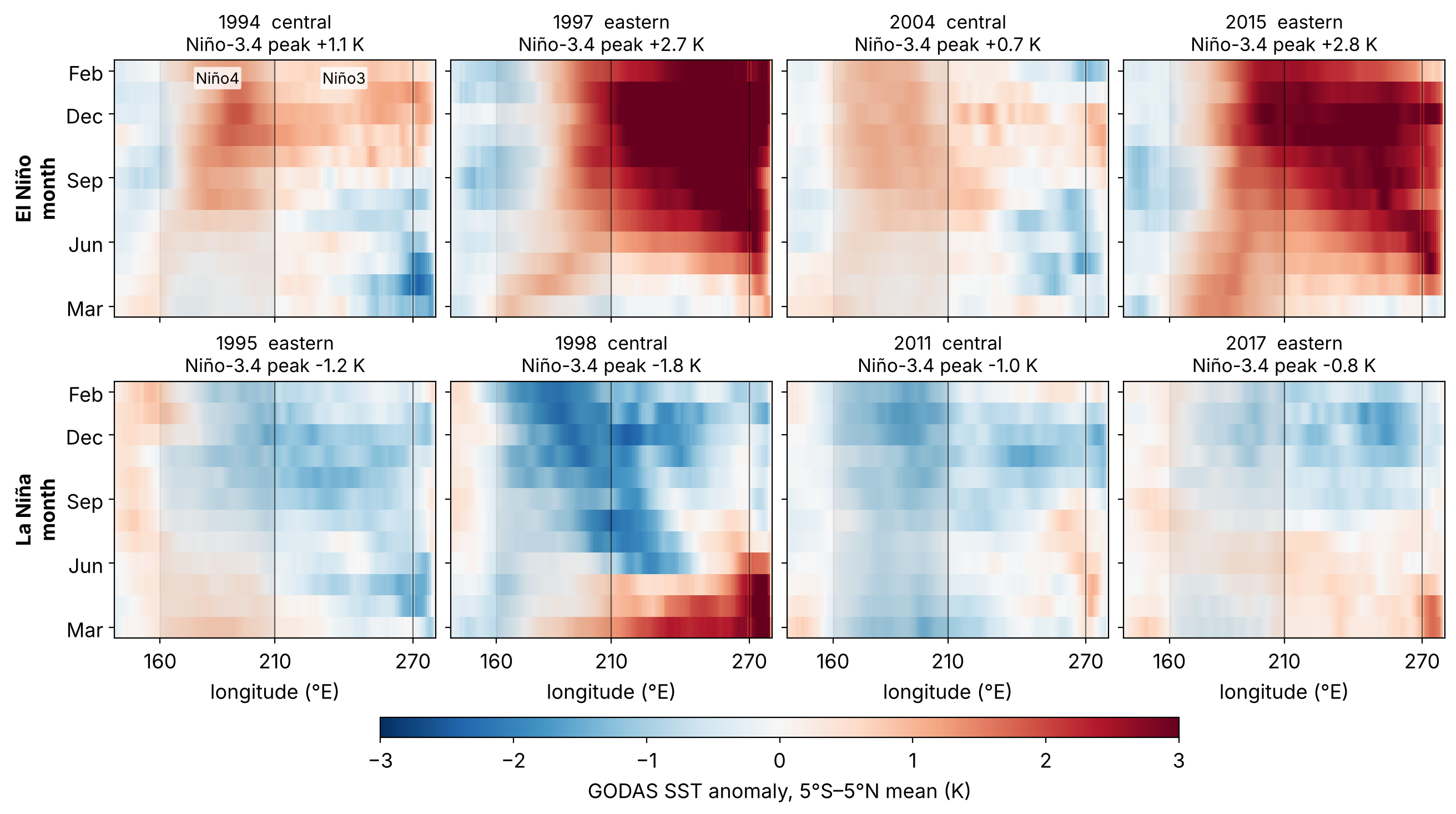}
\caption{GODAS SST anomalies (K), averaged over 5$^\circ$S--5$^\circ$N, for the
eight events of Table~\ref{tab:godas_events}. Top row: El~Ni\~no; bottom row: La~Ni\~na. Each panel shows the 12-month fingerprint window from March of the peak year to the following February. The shaded band marks the Ni\~no-4 region (160--210$^\circ$E) and the solid lines bound the Ni\~no-3 region (210--270$^\circ$E). Colors saturate at $\pm3$\,K.}
\label{fig:godas_hov}
\end{figure}

Figure~\ref{fig:godas_hov} shows the GODAS SST anomalies for the eight events
of Table~\ref{tab:godas_events}, averaged over 5$^\circ$S--5$^\circ$N and plotted over the 12-month window from March of the peak year used for the fingerprints. Anomalies are relative to the 1981--2020 monthly climatology of the same record. The eastern events (1997 and 2015 El~Ni\~no, 1995 and 2017 La~Ni\~na) show their largest anomalies east of 210$^\circ$E, in the Ni\~no-3 region, while the central events (1994 and 2004 El~Ni\~no, 1998 and 2011 La~Ni\~na) show theirs between 160$^\circ$E and 210$^\circ$E, in the Ni\~no-4 region. The 1994 and 2004 El~Ni\~no events are weak, and the 1994 event has the smallest flavor margin of the set (Table~\ref{tab:godas_events}), so its pattern is the least distinct. The windows are fixed in length, so they can include parts of a neighbouring event: the first months of the 1998 La~Ni\~na window still show the end of the 1997--98 El~Ni\~no. Hence a fixed-window fingerprint need not isolate a single event's mechanism-related signals. Only GODAS SST was assimilated in this application; independent subsurface and current observations were not used to validate the reconstructed mechanisms. 

The SST patterns do not independently verify thermocline, wind, or advective dominance in individual events. Such interpretations would require additional physical diagnostics. The present LLM comparison measures agreement with the inferred fingerprint.

\end{document}